\documentclass{gtart_h}
\def\ifplaintex{\expandafter\ifx\csname documentclass\endcsname\relax}

\def\gtp{{\mathsurround=0pt\it $\cal G\mskip-2mu$eometry \&\ 
$\cal T\!\!$opology $\cal P\!$ublications}}  % GT publications

\def\recd{{\small Received:\qua\receiveddate\ifx\reviseddate\relax
\else\qquad Revised:\qua\reviseddate\fi\par}} 

\def\lognumber#1{\def\thelognumber{#1}}
\def\volumenumber#1{\def\thevolumenumber{#1}}
\def\volumeyear#1{\def\thevolumeyear{#1}}
\def\papernumber#1{\def\thepapernumber{#1}}
\def\pagenumbers#1#2{\def\startpage{#1}\def\finishpage{#2}}
\def\published#1{\def\publishdate{#1}}

\def\received#1{\def\receiveddate{#1}}
\def\revised#1{\def\reviseddate{#1}}
\def\accepted#1{\def\accepteddate{#1}}

\def\asciiemail#1{\def\theasciiemail{#1}}

\long\def\asciiabstract#1{\long\def\theasciiabstract{#1}}
\def\asciikeywords#1{\def\theasciikeywords{#1}}

\let\\\par\let\thelognumber\relax\let\thevolumenumber\relax
\let\thepapernumber\relax\let\thevolumeyear\relax\let\startpage\relax
\let\finishpage\relax\let\publishdate\relax\let\receiveddate\relax
\let\reviseddate\relax\let\accepteddate\relax\let\theasciititle\relax
\let\theasciiauthors\relax
\let\theasciiabstract\relax\let\theasciikeywords\relax

\let\theasciiemail\relax

\ifplaintex
\font\logobig=cmssbx10 scaled 3836
\font\logomed=cmssbx10 scaled 2557
\else
\font\logobig=cmssbx10 scaled 4200
\font\logomed=cmssbx10 scaled 2800
\fi

\long\def\makeagttitle{   %%% start of definition of \makeagttitle
\count0=\startpage
\agt\hfill      %   Journal title (top left) 
\hbox to 45truept{\vbox to 0pt{\vglue -13truept{\logomed A\kern -.37em{\logobig 
T}\kern -.38em G}\vss}\hss}
\break
{\small Volume \thevolumenumber\ (\thevolumeyear)
\startpage--\finishpage\nl
Published: \publishdate}

\vglue .25truein

{\parskip=0pt\leftskip 0pt plus
1fil\def\\{\par\smallskip}{\Large\bf\thetitle}\par\medskip} \vglue
0.05truein

{\parskip=0pt\leftskip 0pt plus 1fil\def\\{\par}{\sc\theauthors}
\par\medskip}%
 
\vglue 0.03truein 

{\small\leftskip 25truept\rightskip 25truept{\bf Abstract}\stdspace\theabstract

{\bf AMS Classification}\stdspace\theprimaryclass
\ifx\thesecondaryclass\relax\else; \thesecondaryclass\fi\par
{\bf Keywords}\stdspace \thekeywords\par}\vglue 7truept

}   %%%% end of definition of \makeagttitle

\ifplaintex
\font\phead=cmsl9 scaled 950
\font\pnum=cmbx10 scaled 913
\font\pfoot=cmsl9 scaled 950
\headline{\vbox to 0pt{\vskip -4.5mm\line{\small\phead\ifnum
\count0=\startpage ISSN 1472-2739 (on-line) 1472-2747 (printed)
\hfill {\pnum\folio}\else\ifodd\count0\def\\{ }% 
\ifx\theshorttitle\relax\thetitle\else\theshorttitle\fi\hfill{\pnum\folio}
\else\def\\{ and }{\pnum\folio}\hfill\ifx\theshortauthors\relax\theauthors
\else\theshortauthors\fi\fi\fi}\vss}}
\footline{\vbox to 0pt{\vglue 0mm\line{\small\pfoot\ifnum\count0=\startpage
\copyright\ \gtp\hfill\else
\agt, Volume \thevolumenumber\ (\thevolumeyear)\hfill\fi}\vss}}
\else
\font=cmsl9 scaled 1050
\font\lnum=cmbx10 
\font=cmsl9 scaled 1050
\makeatletter
\def\@oddhead{{\small\lhead\ifnum\count0=\startpage ISSN 1472-2739 
(on-line) 1472-2747 (printed)\hfill {\lnum\number\count0}\else\ifodd\count0
\def\\{ }\ifx\theshorttitle\relax \thetitle \else\theshorttitle\fi\hfill
{\lnum\number\count0}\else\def\\{ and }{\lnum\number\count0}
\hfill\ifx\theshortauthors\relax 
\theauthors\else\theshortauthors\fi\fi\fi}}\def\@evenhead{\@oddhead}
\def\@oddfoot{\small\lfoot\ifnum\count0=\startpage\copyright\ \gtp\hfill\else
\agt, Volume \thevolumenumber\ (\thevolumeyear)\hfill\fi}
\def\@evenfoot{\@oddfoot}
\makeatother
\fi
\let\maketitlepage\makeagttitle
\let\makeshorttitle\maketitlepage
\let\maketitle\maketitlepage

\newwrite\gtoutfile
\long\gdef\makeheadfile{  %%% start of definition of \makeheadfile
{\def\\{, }\def\s{ }
\immediate\openout\gtoutfile head.xxx
\immediate\write\gtoutfile{Proxy-for: \ifx\theasciiauthors\relax
\theauthors\else\theasciiauthors\fi\s<\ifx\theasciiemail\relax\theemail\else\theasciiemail\fi>}
\immediate\write\gtoutfile{\noexpand\\}
\immediate\write\gtoutfile{Authors: \ifx\theasciiauthors\relax
\theauthors\else\theasciiauthors\fi}
{\def\\{ }\immediate\write\gtoutfile{Title: \ifx\theasciititle\relax
\thetitle\else\theasciititle\fi}}
\immediate\write\gtoutfile{Subj-class: GT or SG, GR etc}
\immediate\write\gtoutfile{MSC-class: \theprimaryclass\ifx\thesecondaryclass\relax\else, \thesecondaryclass\fi}
\immediate\write\gtoutfile{Journal-ref: Algebr. Geom. Topol. \thevolumenumber\s
(\thevolumeyear) \startpage-\finishpage}
\immediate\write\gtoutfile{Comments: Published by Algebraic and
Geometric Topology at}
\immediate\write\gtoutfile{\s\s\s  http://www.maths.warwick.ac.uk/agt/AGTVol\thevolumenumber/agt-\thevolumenumber-\thepapernumber.abs.html}
\immediate\write\gtoutfile{\noexpand\\}
\immediate\write\gtoutfile{}
\ifx\theasciiabstract\relax
\immediate\write\gtoutfile{\theabstract}\else
\immediate\write\gtoutfile{\theasciiabstract}\fi
\immediate\write\gtoutfile{}
\immediate\write\gtoutfile{\noexpand\\}
\immediate\write\gtoutfile{}
\immediate\closeout\gtoutfile}}  %%% end of definition of \makeheadfile

\def\maketitlepage{\makeagttitle\makeheadfile}
\let\makeshorttitle\maketitlepage
\let\maketitle\maketitlepage

\lognumber{27}
\volumenumber{5}
\volumeyear{2005}
\papernumber{27}
\pagenumbers{653}{690}
\received{2 June 2004} 
\revised{3 June 2005}
\accepted{21 June 2005}
\published{5 July 2005}

\usepackage{amssymb} 
\usepackage{amsmath}%,amscd} 

\input xy
\xyoption{all}

\SelectTips{cm}{}
\newdir{ >}{{}*!/-6pt/@{>}}       

\newcommand{\C}{\mathbb C}
\newcommand{\Ext}{\operatorname{Ext}}
\newcommand{\F}{\mathbb F}
\newcommand{\R}{\mathbb R}
\newcommand{\T}{\mathbb T}
\newcommand{\Z}{\mathbb Z}
\newcommand{\colim}{\operatornamewithlimits{colim}}

\newcommand{\tmf}{tm\!f}
\renewcommand{\phi}{\varphi}
\newcommand{\THH}{\hbox{\it THH\/}}

\numberwithin{equation}{section}

\newtheorem{thm}{Theorem}[section]
\newtheorem{lemma}[thm]{Lemma}
\newtheorem{prop}[thm]{Proposition}
\newtheorem{cor}[thm]{Corollary}

\theoremstyle{definition}
\newtheorem{remark}[thm]{Remark}   

\begin{document}

\title{Differentials in the homological\\homotopy fixed point spectral sequence}
\authors{Robert R. Bruner\\John Rognes}
\address{Department of Mathematics, Wayne State University, 
Detroit, MI 48202, USA\\Department of Mathematics, University of Oslo,
NO-0316 Oslo, Norway}
\gtemail{\mailto{rrb@math.wayne.edu}, \mailto{rognes@math.uio.no}}
\asciiemail{rrb@math.wayne.edu, rognes@math.uio.no}

\begin{abstract}   % type your abstract below
We analyze in homological terms the homotopy fixed point spectrum of a
$\T$-equivariant commutative $S$-algebra $R$.  There is a homological
homotopy fixed point spectral sequence with $E^2_{s,t} =
H^{-s}_{gp}(\T; H_t(R; \F_p))$, converging conditionally to the
continuous homology $H^c_{s+t}(R^{h\T}; \F_p)$ of the homotopy fixed
point spectrum.  We show that there are Dyer--Lashof operations
$\beta^\epsilon Q^i$ acting on this algebra spectral sequence, and
that its differentials are completely determined by those originating
on the vertical axis.  More surprisingly, we show that for each
class~$x$ in the $E^{2r}$-term of the spectral sequence there are $2r$
other classes in the $E^{2r}$-term (obtained mostly by Dyer--Lashof
operations on $x$) that are infinite cycles, i.e., survive to the
$E^\infty$-term.  We apply this to completely determine the
differentials in the homological homotopy fixed point spectral
sequences for the topological Hochschild homology spectra $R = \THH(B)$
of many $S$-algebras, including $B = MU$, $BP$, $ku$, $ko$ and~$\tmf$.
Similar results apply for all finite subgroups $C \subset \T$, and for
the Tate- and homotopy orbit spectral sequences.  This work is part of
a homological approach to calculating topological cyclic homology and
algebraic $K$-theory of commutative $S$-algebras.
\end{abstract}

\asciiabstract{%
We analyze in homological terms the homotopy fixed point spectrum of a
T-equivariant commutative S-algebra R.  There is a homological
homotopy fixed point spectral sequence with E^2_{s,t} = H^{-s}_{gp}(T;
H_t(R; F_p)), converging conditionally to the continuous homology
H^c_{s+t}(R^{hT}; F_p) of the homotopy fixed point spectrum.  We show
that there are Dyer-Lashof operations beta^epsilon Q^i acting on this
algebra spectral sequence, and that its differentials are completely
determined by those originating on the vertical axis.  More
surprisingly, we show that for each class x in the $^{2r}-term of the
spectral sequence there are 2r other classes in the E^{2r}-term
(obtained mostly by Dyer-Lashof operations on x) that are infinite
cycles, i.e., survive to the E^infty-term.  We apply this to
completely determine the differentials in the homological homotopy
fixed point spectral sequences for the topological Hochschild homology
spectra R = THH(B) of many S-algebras, including B = MU, BP, ku, ko
and tmf.  Similar results apply for all finite subgroups C of T, and
for the Tate- and homotopy orbit spectral sequences.  This work is
part of a homological approach to calculating topological cyclic
homology and algebraic K-theory of commutative S-algebras.}



\asciikeywords{Homotopy fixed points, Tate spectrum, homotopy orbits, 
commutative S-algebra, Dyer-Lashof operations, differentials,
topological Hochschild homology, topological cyclic homology,
algebraic K-theory}

\maketitle 

\section{Introduction}

In this paper we study a homological version of the homotopy fixed
point spectral sequence of an equivariant spectrum, with emphasis on the
interaction between differentials in the spectral sequence and strictly
commutative products in the equivariant spectrum.  We begin by describing
the intended application of this study.

Consider a connective $S$-algebra $B$, such as the sphere spectrum $S$,
the complex cobordism spectrum $MU$ or the Eilenberg--Mac\,Lane spectrum
$H\Z$ of the integers, in either of the current frameworks \cite{BHM93},
\cite{EKMM97}, \cite{HSS00}, \cite{MMSS01} for structured ring spectra.
The algebraic $K$-theory spectrum $K(B)$ can, according to the main
theorem of \cite{Du97}, be very well approximated by the topological
cyclic homology spectrum $TC(B)$ of \cite{BHM93}.  The latter is
constructed from the $\T$-equivariant topological Hochschild homology
spectrum $X = \THH(B)$, where $\T$ is the circle group, as a homotopy
limit of the fixed point spectra $X^C$, suitably indexed over finite
cyclic subgroups $C \subset \T$.

These fixed point spectra are in turn often well approximated, via
canonical maps $\gamma \colon X^C \to X^{hC}$, by homotopy fixed point
spectra $X^{hC} = F(EC_+, X)^C$.  In principle, the homotopy groups of
the latter can be computed by the {\it homotopical\/} homotopy fixed
point spectral sequence
\begin{equation}
E^2_{s,t} = H^{-s}_{gp}(C; \pi_t(X)) \Longrightarrow \pi_{s+t}(X^{hC}) \,,
\label{e1.1}
\end{equation}
which is derived by applying homotopy groups to the tower of fibrations,
with limit~$X^{hC}$, that arises from the equivariant skeleton filtration
on the free contractible $C$-space $EC$.

Such computations presume a detailed knowledge of the homotopy groups
$\pi_*(X)$ of the $\T$-equivariant spectrum in question.  For example,
the papers \cite{HM03} and~\cite{AuR02} deal with the cases when $B$
is the valuation ring $H{\mathcal{O}}_K$ in a local number field $K$ and the
Adams summand $\ell_p$ in $p$-complete connective topological $K$-theory
$ku_p$, respectively.  In most other cases the spectral sequence~(\ref{e1.1})
cannot be calculated, because the homotopy groups $\pi_*(X)$ are not
sufficiently well known.

It happens more frequently in stable homotopy theory that we are familiar
with the homology groups $H_*(X; \F_p)$.  Applying mod~$p$ homology,
rather than homotopy, to the tower of fibrations that approximates
$X^{hC}$ leads to a {\it homological\/} homotopy fixed point spectral
sequence
\begin{equation}
E^2_{s,t} = H^{-s}_{gp}(C; H_t(X; \F_p)) \Longrightarrow H^c_{s+t}(X^{hC};
\F_p) \,.
\label{e1.2}
\end{equation}
This spectral sequence converges conditionally \cite{Bo99} to the
(inverse) limit of the resulting tower in homology, which is not
$H_*(X^{hC}; \F_p)$, but a ``continuous'' version $H^c_*(X^{hC}; \F_p)$
of it, for homology does not usually commute with limits.

This continuous homology, considered as a comodule over the dual
Steenrod algebra $A_*$ \cite{Mi58}, is nonetheless a powerful invariant
of $X^{hC}$.  In particular, when $X$ is bounded below and of finite
type there is a strongly convergent spectral sequence
\begin{equation}
E_2^{s,t} = \Ext_{A_*}^{s,t}(\F_p, H^c_*(X^{hC}; \F_p))
\Longrightarrow \pi_{t-s}(X^{hC})^\wedge_p
\label{e1.3}
\end{equation}
which can be obtained as an inverse limit of Adams spectral sequences
\cite[7.1]{CMP87}.  Hence the continuous homology does in some sense
determine the $p$-adic homotopy type of $X^{hC}$.

A form of the spectral sequence~(\ref{e1.3}) was most notably applied in the
proofs by W.~H.~Lin \cite{LDMA80} and J.~Gunawardena \cite{AGM85} of
the Segal conjecture for cyclic groups of prime order.  The conjecture
corresponds to the special case of the discussion above when $B = S$
is the sphere spectrum, so $X = \THH(S) = S$ is the $\T$-equivariant
sphere spectrum, which is split \cite[II.8]{LMS86} so that $X^{hC}
\simeq F(BC_+, S) = D(BC_+)$.  The proven Segal conjecture \cite{Ca84}
then tells us that for each $p$-group $C$ the comparison map $\gamma \colon
X^C \to X^{hC}$ is a $p$-adic equivalence.  The proof of the general
(cyclic) case is by reduction to the initial case when $C = C_p$
is of prime order, and therefore relies on the theorems of Lin and
Gunawardena cited above.  In this case, of course, we do not know
$\pi_*(X) = \pi_*(S)$ sufficiently well, but $H_*(X; \F_p) = \F_p$ is
particularly simple.  The proof of the theorems of Lin and Gunawardena
now amounts to showing that although the natural homomorphism $\gamma_*
\colon H_*(X^C; \F_p) \to H^c_*(X^{hC}; \F_p)$ of $A_*$-comodules is not in
itself an isomorphism, it does induce an isomorphism of $E_2$-terms upon
applying the functor $\Ext_{A_*}^{**}(\F_p, -)$.

Returning to the general situation, we are therefore interested in
studying (i)~the differentials in the homological homotopy fixed point
spectral sequence~(\ref{e1.2}) above, and (ii)~the $A_*$-comodule extension
questions relating its $E^\infty$-term to the abutment $H^c_*(X^{hC};
\F_p)$.  There will in general be non-trivial differentials in~(\ref{e1.2}),
but our main Theorem~\ref{t1.5} below provides a very general and useful
vanishing result, as is illustrated by the examples in Section~6.
The identification of the $A_*$-comodule structure on the abutment
plays a crucial role already in the case $X = S$, but requires further
study beyond that given here, and will be presented in the forthcoming
Ph.D. thesis of Sverre Lun{\o}e--Nielsen \cite{L-N}.

Now suppose that $B$ is a commutative $S$-algebra, in either one of
the structured categories listed at the outset.  Then $X = \THH(B)$ can
be constructed as a $\T$-equivariant commutative $S$-algebra, which we
hereafter denote $R = X$, in the naive sense of a commutative $S$-algebra
with a continuous $\T$-action through commutative $S$-algebra maps.
Starting at this point we need to be in a technical framework where
naively $\T$-equivariant commutative $S$-algebras $R$ make sense, and
either of the $S$-modules of \cite{EKMM97}, the orthogonal spectra of
\cite{MMSS01}, the equivariant orthogonal spectra of \cite{MM02} or the
topological version of the symmetric spectra of \cite{HSS00} will do.
To be concrete, we can follow \cite[Ch.~IX]{EKMM97}.

Then the tower of fibrations with limit $R^{hC} = X^{hC}$ is one
of commutative $S$-algebras.  It follows that there are Dyer--Lashof
operations acting on the spectral sequence~(\ref{e1.2}) in this case, somewhat
analogously to the action by Steenrod operations in the Adams spectral
sequence of a commutative $S$-algebra \cite[Ch.~IV]{BMMS86}.  In the
latter case there are very interesting relations between the Adams
differentials and the Steenrod operations, which propagate early
differentials to higher degrees.  The initial motivation for the
present article was to determine the analogous interaction between the
differentials and the Dyer--Lashof operations in the homological homotopy
fixed point spectral sequence of a commutative $S$-algebra.  However,
the analogy with the behavior of differentials in the Adams spectral
sequence is more apparent than real, suggesting neither the survival
to $E^\infty$ of some classes, nor the method of proof.  In particular,
there is no analog in the Adams spectral sequence of our main vanishing
result, Theorem~\ref{t1.5}.

For each finite subgroup $C \subset \T$ the homological spectral sequence
for $R^{hC}$ is an algebra over the corresponding homological spectral
sequence for $R^{h\T}$, as outlined in Section~7, so it will suffice
for us to consider the circle homotopy fixed points $R^{h\T}$ and the
case $C = \T$ of the spectral sequence~(\ref{e1.2}).  Our first results in
Sections~2--4 can then be summarized as follows.

\begin{thm}
\label{t1.4}
{\rm(a)}\qua
Let $R$ be a $\T$-equivariant commutative $S$-algebra.  Then there is
a natural $A_*$-comodule algebra spectral sequence
$$
E^2_{**}(R) = H^{-*}_{gp}(\T; H_*(R; \F_p)) = P(y) \otimes H_*(R; \F_p)
$$
with $y$ in bidegree $(-2,0)$, converging conditionally to the continuous
homology
$$
H^c_*(R^{h\T}; \F_p) = \lim_n H_*(F(S(\C^n)_+, R)^{\T}; \F_p)
$$
of the homotopy fixed point spectrum $R^{h\T} = F(E\T_+, R)^{\T}$.

{\rm(b)}\qua
There are natural Dyer--Lashof operations $\beta^\epsilon Q^i$ acting
vertically on this homological homotopy fixed point spectral sequence.
For each element $x \in E^{2r}_{0,t} \subset H_t(R; \F_p)$ we have
the relation
$$
d^{2r}(\beta^\epsilon Q^i(x)) = \beta^\epsilon Q^i(d^{2r}(x))
$$
for every integer $i$ and $\epsilon \in \{0,1\}$.  If $d^{2r}(x) = y^r
\cdot \delta x$ with $\delta x \in H_{t+2r-1}(R; \F_p)$, then the right
hand side $\beta^\epsilon Q^i(d^{2r}(x))$ is defined to be $y^r \cdot
\beta^\epsilon Q^i(\delta x)$.

{\rm(c)}\qua
The classes $y^n$ are infinite cycles, so the differentials from the
vertical axis $E^{2r}_{0,*}$ propagate to each column by the relation
$$
d^{2r}(y^n \cdot x) = y^n \cdot d^{2r}(x)
$$
for every $x \in E^{2r}_{0,*}$, $2r\ge2$, $n\ge0$.  Hence there are
isomorphisms $E^{2r}_{**} \equiv P(y) \otimes E^{2r}_{0,*}$ for all
$2r\ge2$, modulo possible $y$-torsion in filtrations $-2r+4 \le s \le 0$.
\end{thm}

For proofs, see Propositions~\ref{p2.6}, \ref{p4.1} and~\ref{p4.4}, and Lemma~\ref{l4.6}.  The key
idea is to identify the differentials in the homological homotopy fixed
point spectral sequence as obstructions to extending equivariant maps, as
explained in Section~3.  Note that the spectral sequence is concentrated
in even columns, so all differentials of odd length ($d^r$ with $r$ odd)
must vanish.

Our main theorem is the following vanishing result for differentials.

\begin{thm}
\label{t1.5}
Let $R$ be a $\T$-equivariant commutative $S$-algebra, suppose that
$x \in H_t(R; \F_p)$ survives to the $E^{2r}$-term $E^{2r}_{0,t}
\subset H_t(R; \F_p)$ of the homological homotopy fixed point spectral
sequence for $R$, and write $d^{2r}(x) = y^r \cdot \delta x$.

{\rm(a)}\qua
For $p=2$, the $2r$ classes
$$
x^2 = Q^t(x),\  Q^{t+1}(x),\  \dots \ ,\  Q^{t+2r-2}(x)
\quad\text{and}\quad Q^{t+2r-1}(x) + x \delta x
$$
all survive to the $E^\infty$-term, i.e., are infinite cycles.

{\rm(b)}\qua
For $p$ odd and $t = 2m$ even, the $2r$ classes
$$
x^p = Q^m(x),\  \beta Q^{m+1}(x),\  \dots \ ,\  Q^{m+r-1}(x)
\quad\text{and}\quad x^{p-1} \delta x
$$
all survive to the $E^\infty$-term, i.e., are infinite cycles.

{\rm(c)}\qua
For $p$ odd and $t = 2m-1$ odd, the $2r$ classes
$$
\beta Q^m(x),\  Q^m(x),\  \dots \ ,\  \beta Q^{m+r-1}(x)
\quad\text{and}\quad Q^{m+r-1}(x) - x(\delta x)^{p-1}
$$
all survive to the $E^\infty$-term, i.e., are infinite cycles.
\end{thm}

This will be proved in Section~5 as Theorem~\ref{t5.1}.  There are similar
extensions of our results to the Tate constructions $R^{tC} =
[\widetilde{EC} \wedge F(EC_+, R)]^C$ and homotopy orbit spectra $R_{hC}
= EC_+ \wedge_C R$, but to keep the exposition simple these are also
only discussed in Section~7.

\begin{remark}
In more detail, in the three cases of the theorem the infinite cycles
can be listed as follows:
\begin{enumerate}
\item[(a)]
$x^2 = Q^t(x)$, $Q^i(x)$ for $t+1 \le i \le t+2r-2$ and $Q^{t+2r-1}(x)
+ x \delta x$;
\item[(b)]
$x^p = Q^m(x)$, $\beta^\epsilon Q^i(x)$ for $m+1 \le i \le m+r-1$,
$\epsilon \in \{0, 1\}$ and $x^{p-1} \delta x$;
\item[(c)]
$\beta^\epsilon Q^i(x)$ for $m \le i \le m+r-2$, $\epsilon \in
\{0, 1\}$, $\beta Q^{m+r-1}(x)$ and $Q^{m+r-1}(x) - x(\delta x)^{p-1}$.
\end{enumerate}

In particular, for any element $x \in H_t(R; \F_p)$, with $d^2(x) =
y \cdot \delta x$, there are two infinite cycles
\begin{enumerate}
\item[(a)]
$x^2 = Q^t(x)$ and $Q^{t+1}(x) + x \delta x$ when $p=2$,
\item[(b)]
$x^p = Q^m(x)$ and $x^{p-1} \delta x$ when $p$ is odd and $t=2m$, and
\item[(c)]
$\beta Q^m(x)$ and $Q^m(x) - x(\delta x)^{p-1}$ when $p$ is odd
and $t=2m-1$.
\end{enumerate}
\end{remark}

As applications of our main results, we turn in Section~6 to the
study of the algebraic $K$-theory spectrum $K(MU)$, interpolating
between $K(S)$ (which is Waldhausen's $A(*)$, related to high dimensional
geometric topology) and $K(\Z)$ (which relates to the Vandiver and Leopoldt
conjectures, and other number theory), by the methods of topological
cyclic homology.  Hence we must study the fixed- and homotopy fixed
point spectra of the commutative $S$-algebra $R = \THH(MU)$, for various
subgroups $C$ of the circle group $\T$.  It is known that $H_*(MU;
\F_p) = P(b_k \mid k\ge1)$, where $P(-)$ denotes the polynomial algebra
over $\F_p$ and $|b_k| = 2k$, from which it follows (\cite[4.3]{MS93}
or \cite{CS}) that $H_*(\THH(MU); \F_p) = H_*(MU; \F_p) \otimes E(\sigma
b_k \mid k\ge1)$, where $E(-)$ denotes the exterior algebra over $\F_p$
and $\sigma \colon H_*(R; \F_p) \to H_{*+1}(R; \F_p)$ is the degree~$+1$
operator induced by the circle action.  Hence the homological homotopy
fixed point spectral sequence for $\THH(MU)^{h\T}$ begins
$$
E^2_{**} = P(y) \otimes P(b_k \mid k\ge1) \otimes E(\sigma b_k \mid k\ge1)
\,.
$$
There are differentials $d^2(b_k) = y \cdot \sigma b_k$ for all $k\ge1$,
so by our Theorem~\ref{t1.4}
$$
E^4_{**} = P(y) \otimes P(b_k^p \mid k\ge1) \otimes E(b_k^{p-1}
\sigma b_k \mid k\ge1)
$$
plus some classes (the image of $\sigma$) in filtration $s=0$.
By our Theorem~\ref{t1.5} the spectral sequence collapses completely at
the $E^4$-term, so that
$$
H^c_*(\THH(MU)^{h\T}; \F_p) = P(y) \otimes P(b_k^p \mid k\ge1)
\otimes E(b_k^{p-1} \sigma b_k \mid k\ge1)
$$
plus some classes in filtration zero, as an algebra.  The identification
of the $A_*$-comodule extensions remains, for which we refer to
the cited Ph.D. thesis \cite{L-N}.  This provides the input for
the inverse limit of Adams spectral sequences~(\ref{e1.3}) converging to
$\pi_*(\THH(MU)^{h\T})^\wedge_p$, which approximates the topological
version $TF(MU)$ of negative cyclic homology, and which determines the
topological cyclic homology of $MU$ by a fiber sequence
\[
\xymatrix{
TC(MU) \ar[r]^{\pi} & TF(MU) \ar[r]^{R-1} & TF(MU) \,.
}
\]
The fiber of the cyclotomic trace map $K(MU) \to TC(MU)$ is equivalent
to that of $K(\Z) \to TC(\Z)$, by \cite{Du97}, which now is quite well
known \cite{Ro02}, \cite{Ro03}.  Our theorem therefore provides a key
input to the computation of $K(MU)$.  See Theorem~\ref{t6.4}(a).

Similar applications are given for the connective Johnson--Wilson spectra
$B = BP\langle n\rangle$, for $p$ and $n$ such that these are commutative
$S$-algebras, and the (higher real) commutative $S$-algebras $B = ko$ and
$\tmf$ for $p=2$.  See Section~6.  Lastly, we can also show the collapse
at the $E^4$-term of the homological homotopy fixed point spectral
sequence for $R = \THH(BP)$, where $BP$ is the $p$-local Brown--Peterson
$S$-algebra \cite{BJ02}, without making the (presently uncertain)
assumption that $BP$ can be realized as a commutative $S$-algebra.
See Theorem~\ref{t6.4}(b).  This is possible by the homological approach, since
the split surjection $H_*(MU; \F_p) \to H_*(BP; \F_p)$ prevails throughout
the homological spectral sequences.

\section{ A homological spectral sequence }

In this section we construct the homological homotopy fixed point spectral
sequence for a spectrum with a circle action.

Let $\T \subset \C^*$ be the circle group.  As our concrete model
for a free contractible $\T$-CW complex $E\T$ we take the unit sphere
$S(\C^\infty) \subset \C^\infty$ with the usual coordinatewise action
by $\T$.  It has one $\T$-equivariant cell in each even non-negative
dimension $2n \ge 0$.  The equivariant $2n$-skeleton is the odd
$(2n+1)$-sphere $E\T^{(2n)} = S(\C^{n+1})$, which is obtained from the
equivariant $(2n-2)$-skeleton $E\T^{(2n-2)} = S(\C^n)$ by attaching a free
$\T$-equivariant $2n$-cell $\T \times D^{2n}$ along the group action map
\begin{equation}
\alpha \colon \T \times \partial D^{2n} \to S(\C^n) \,.
\label{e2.1}
\end{equation}
This map is $\T$-equivariant when we give $\partial D^{2n}$ the trivial
$\T$-action and $S(\C^n)$ the free $\T$-action.  Hence there is a
$\T$-equivariant filtration
$$
\emptyset \subset S(\C) \subset \dots \subset S(\C^n) \subset S(\C^{n+1})
\subset \dots
$$
with colimit $E\T$, and $\T$-equivariant cofiber sequences
$$
S(\C^n) \to S(\C^{n+1}) \to \T_+ \wedge S^{2n}
$$
for $n\ge0$.  Here $\T$ acts trivially on $S^{2n} = D^{2n}/\partial
D^{2n}$.

Let $X$ be any spectrum with $\T$-action, i.e., a naively $\T$-equivariant
spectrum.  The {\sl homotopy fixed point spectrum} of $X$ is defined as
the mapping spectrum
$$
X^{h\T} = F(E\T_+, X)^{\T}
$$
of $\T$-equivariant based maps from $E\T_+$ to $X$.  Here $E\T_+$
should be interpreted as the unreduced suspension spectrum of $E\T$.
The $\T$-equivariant filtration of $E\T = S(\C^\infty)$ induces a tower
of fibrations
\begin{equation}
\dots \to F(S(\C^{n+1})_+, X)^{\T} \to F(S(\C^n)_+, X)^{\T} \to
\dots \to F(S(\C)_+, X)^{\T} \to *
\label{e2.2}
\end{equation}
with the homotopy fixed point spectrum as its limit
$$
X^{h\T} = \lim_n F(S(\C^n)_+, X)^{\T} \,,
$$
weakly equivalent to the homotopy limit.  The cofiber sequences above
induce cofiber sequences of spectra
$$
\Sigma^{-2n} X = F(\T_+ \wedge S^{2n}, X)^{\T}
\xrightarrow{k} F(S(\C^{n+1})_+, X)^{\T} \xrightarrow{i} F(S(\C^n)_+, X)^{\T}
$$
for each $n\ge0$.  We place $F(E\T^{(-s-1)}_+, X)^{\T}$ in filtration $s$,
for each $s \in \Z$, and obtain a chain of cofiber sequences of spectra:
$$
\xymatrix{
X^{h\T} \ar[r]
	& F(S(\C^{n+1})_+, X)^{\T} \ar[r]^i
	& F(S(\C^n)_+, X)^{\T} \ar[d]^j \ar[r]^i
	& F(S(\C^n)_+, X)^{\T} \ar[d]^j \\
& \dots & \Sigma^{-2n} X \ar[ul]^k & {*} \ar[ul]^k
}
$$
Here the filtrations $s = -2n-1$, $-2n$ and $-2n+1$ are displayed,
and the vertical maps labeled $j$ are of degree~$-1$.

Now apply mod~$p$ homology $H_*(-; \F_p)$ to this chain of cofiber
sequences, to obtain a homologically indexed unrolled exact couple
\cite[0.1]{Bo99} of graded $A_*$-comodules
$$
\xymatrix{
A_{-\infty} \ar[r] & A_{-2n-1} \ar[r]^i
	& A_{-2n} \ar[r]^i \ar[d]^j & A_{-2n+1} \ar[r] \ar[d]^j
	& A_\infty \\
& \dots & E_{-2n} \ar[ul]^k & E_{-2n+1} \ar[ul]^k 
}
$$
with
$$
A_{s,t} = H_{s+t}(F(E\T^{(-s-1)}_+, X)^{\T}; \F_p)
$$
equal to $H_{s+t}(F(S(\C^n)_+, X)^{\T}; \F_p)$ for $s = -2n$ and $s =
-2n+1$, and
$$
E_{s,t} = H_{s+t}(F(E\T^{(-s)}/E\T^{(-s-1)}, X)^{\T}; \F_p)
$$
equal to $H_{s+t}(\Sigma^{-2n} X; \F_p) \cong H_t(X; \F_p)$ for $s =
-2n \le 0$, and zero otherwise.  We write $y^n \cdot x \in E_{-2n,t}$ for
the class that corresponds to $x \in H_t(X; \F_p)$ under this canonical
suspension isomorphism.

The associated spectral sequence $E^r = E^r(X)$ with
$$
E^1_{s,t} = E_{s,t} \cong \begin{cases}
H_t(X; \F_p) & \text{for $s = -2n \le 0$} \\
0 & \text{otherwise,}
\end{cases}
$$
is by definition the homological homotopy fixed point spectral sequence
for the $\T$-equivariant spectrum~$X$.  The $E^1$-term is concentrated
in the non-positive even columns, so $d^1 = 0$ and $E^1 = E^2$.

Alternatively, there is (the finite part of) a Cartan--Eilenberg system
\cite[XV.7]{CE56}, with graded groups
$$
H(a,b) = H_*(F(E\T^{(b-1)}/E\T^{(a-1)}, X)^{\T}; \F_p)
$$
for all integers $a \le b$, which via $E^1_s = H(-s,-s+1)$ gives rise
to the same spectral sequence as the unrolled exact couple above.
Also define $H(a, \infty) = \lim_b H(a, b)$, so that $H(0, \infty) =
H^c_*(X^{h\T}; \F_p)$.  We shall refer to this formalism when discussing
products in Section~4.

Since $A_s = 0$ for $s \ge 0$ we trivially have $A_\infty = \colim_s
A_s = 0$.  Therefore the associated spectral sequence is conditionally
convergent, by \cite[5.10]{Bo99}, in this case to the limit $A_{-\infty}
= \lim_s A_s$.  Indexing the limit system by $n$ in place of~$s$, the
abutment can be written as
\begin{equation}
H^c_*(X^{h\T}; \F_p) = \lim_n H_*(F(S(\C^n)_+, X)^{\T}; \F_p) \,,
\label{e2.3}
\end{equation}
which we call the {\sl continuous homology} of $X^{h\T}$.  Since homology
rarely commutes with the formation of limits, the canonical map
$$
H_*(X^{h\T}; \F_p) \to H^c_*(X^{h\T}; \F_p)
$$
is usually not an isomorphism.  The Segal conjecture provides striking
examples of this phenomenon.

The spectral sequence will be strongly convergent to $H^c_*(X^{h\T};
\F_p)$ if the criterion $RE^\infty_{**} = 0$ of \cite[7.4]{Bo99} is
satisfied, for which it certainly suffices that in each bidegree $(s,t)$
we have $E^r_{s,t} = E^\infty_{s,t}$ for some finite $r = r(s,t)$.
We recall that in this setting, with convergence to the limit, strong
convergence means that the kernels $F_s = \ker(A_{-\infty} \to A_s)$
form an exhaustive complete Hausdorff filtration of $A_{-\infty} =
H^c_*(X^{h\T}; \F_p)$, and that there are isomorphisms $F_s/F_{s-1} \cong
E^\infty_s$ for all $s$.  When $RE^\infty_{**} = 0$, these isomorphisms
arise from the upper quadrangle (derived from \cite[5.6~and~5.9]{Bo99})
in the following commutative diagram:
\begin{equation}
\xymatrix{
A_{-\infty} \ar[rrr] &&& A_{s-1} \\
F_s \ar@{ >->}[u] \ar@{->>}[r] & F_s/F_{s-1} \ar[r]^-{\cong}
	& E^\infty_s \ar@{ >->}[ur] \\
H(-s,\infty) \ar@{->>}[u] \ar[rr] && Z^\infty_s \ar@{->>}[u] \ar@{ >->}[r]
	& E^1_s \ar[uu]_k
}
\label{e2.4}
\end{equation}
The right hand quadrangle is a pull-back, and the outer rectangle is
obtained by applying $F(-, X)^{\T}$ and continuous homology to the
commutative square below.
$$
\xymatrix{
E\T_+ \ar[d] && E\T^{(-s)}_+ \ar[ll] \ar[d]^k \\
E\T/E\T^{(-s-1)} && E\T^{(-s)}/E\T^{(-s-1)} \ar[ll]
}
$$
When discussing products, we shall find it more convenient to view the
isomorphisms $F_s/F_{s-1} \cong E^\infty_s$ as being determined by the
lower and right hand part of~(\ref{e2.4}).

We now turn to giving a convenient description of the $E^2$-term of
the homological homotopy fixed point spectral sequence.  By the group
cohomology of the circle group $\T$, with coefficients in some discrete
group $M$, we mean the singular cohomology of its classifying space
$B\T$, or equivalently, the $\T$-equivariant cohomology of the universal
space $E\T$:
$$
H^{-*}_{gp}(\T; M) = H^{-*}(B\T; M) = H^{-*}_{\T}(E\T; M) \,.
$$
The latter can be computed from the cellular chain complex given by the
$\T$-equivariant skeletal filtration of $E\T$, cf.~\cite[Ch.~10]{GM95}.
For $M = H_t(X; \F_p) = \pi_t(H\F_p \wedge X)$ we can recognize that
cellular chain complex as the row $(E^1_{*,t}, d^1)$.  This uses the
canonical weak equivalence (of e.g.~\cite[III.1]{LMS86})
\begin{equation}
\nu \colon H\F_p \wedge F(Z, X)^{\T} \to F(Z, H\F_p \wedge X)^{\T}
\label{e2.5}
\end{equation}
for each finite $\T$-CW spectrum $Z$, applied in the cases $Z =
E\T^{(-s)}/E\T^{(-s-1)}$ and $Z = E\T^{(-s+1)}/E\T^{(-s-1)}$ to identify
the groups and the differentials, respectively.  Therefore the $E^2$-term
of the homological spectral can be expressed as
$$
E^2_{s,t} = H^{-s}_{gp}(\T; H_t(X; \F_p))
\cong H^{-s}_{gp}(\T; \F_p) \otimes H_t(X; \F_p) \,.
$$
Furthermore, $H^*_{gp}(\T; \F_p) = P(y)$ is the polynomial algebra on
the Euler class $y \in H^2_{gp}(\T; \F_p)$ of the canonical line bundle
over $B\T$.  Thus
$$
E^2_{**} = P(y) \otimes H_*(X; \F_p)
$$
with $y$ in bidegree $(-2,0)$ and $H_t(X; \F_p)$ in bidegree $(0,t)$.
See \cite[Ch.~14]{GM95} for a discussion of related spectral sequences.

\begin{prop}
\label{p2.6}
There is a natural homological spectral sequence of $A_*$-comodules
$$
E^2_{**}(X) = H^{-*}_{gp}(\T; H_*(X; \F_p)) = P(y) \otimes H_*(X; \F_p)
$$
with $y$ in bidegree $(-2,0)$, converging conditionally to the continuous
homology $H^c_*(X^{h\T}; \F_p)$.  We call this the {\sl homological
homotopy fixed point spectral sequence}.  If $H_*(X; \F_p)$ is finite
in each degree, or the spectral sequence collapses at a finite stage,
then the spectral sequence is strongly convergent.
\end{prop}

\begin{remark}
\label{r2.7}
{\rm(a)}\qua
So far, $E^r(X)$ was just defined as an additive spectral sequence,
but we shall later (in Proposition~\ref{p4.1}) justify the reference to the
algebra structure in $P(y)$ by showing that when $X$ is a $\T$-equivariant
commutative $S$-algebra, then $E^r(X)$ is an algebra spectral sequence.
In the special case when $X = S$ the $E^2$-term is $P(y)$ and the algebra
structure is precisely that of the polynomial algebra $P(y)$.  A general
(naively) $\T$-equivariant spectrum $X$ can be considered as a (naively)
$\T$-equivariant $S$-module \cite[IV.2.8(iv)]{MM02}, and its homological
homotopy fixed point spectral sequence $E^r(X)$ becomes a module spectral
sequence over the collapsing algebra spectral sequence $E^r(S)$ with
$E^2(S) = E^\infty(S) = P(y)$.  In this sense the expression for $E^2(X)$
describes a natural $P(y)$-module structure on $E^r(X)$.

{\rm(b)}\qua
As noted in the introduction, it is rather more traditional to apply
the homotopy group functor $\pi_*(-)$ to the tower of fibrations~(\ref{e2.2}),
to obtain an unrolled exact couple and a conditionally convergent
(homotopical) homotopy fixed point spectral sequence
$$
E^2_{s,t} = H^{-s}_{gp}(\T; \pi_t(X))
\Longrightarrow \pi_{s+t}(X^{h\T}) \,.
$$
However, this is not the spectral sequence that we are considering.
Earlier work by Ch.~Ausoni and the second author \cite[Ch.~4]{AuR02},
as well as recent work by S.~Lun{\o}e-Nielsen (and the second author)
\cite{L-N}, supports the assertion that the homological spectral sequence
is an interesting object.

{\rm(c)}\qua
In view of the natural weak equivalence~(\ref{e2.5}), applied to the various
skeletal filtration quotients $Z = E\T^{(m)}/E\T^{(n)}$, the homological
homotopy fixed point spectral sequence for the $\T$-equivariant spectrum
$X$ is in fact isomorphic to the homotopical homotopy fixed point
spectral sequence for $H\F_p \wedge X$, where $\T$ acts trivially on
the Eilenberg--Mac\,Lane spectrum.  More precisely, they have isomorphic
defining Cartan--Eilenberg systems, in the sense of \cite[XV.7]{CE56}.

{\rm(d)}\qua
In general there can be arbitrarily long differentials in the homological
homotopy fixed point spectral sequence.  For example, when $X = S(\C^r)_+$
with the free $\T$-action, the differentials $d^{2r}_{s,0}$ are nonzero
for all even $s = -2n \le 0$.  The point of Theorem~\ref{t1.5} is that this
rarely happens when $X = R$ is a $\T$-equivariant commutative $S$-algebra.
\end{remark}

\section{ Differentials }

We now make the differentials in the homological homotopy fixed point
spectral sequence more explicit, as obstructions to extending equivariant
maps.

Consider a class $x \in H_t(X; \F_p)$, represented at the $E^2$-term of
the homological spectral sequence in bidegree $(0,t)$.  We briefly write
$H = H\F_p$ for the mod~$p$ Eilenberg--Mac\,Lane spectrum.  Then $x$
can be represented as a non-equivariant map $S^t \to H \wedge X$, or
equivalently as a $\T$-equivariant map
$$
x \colon S(\C)_+ \wedge S^t \to H \wedge X \,.
$$
Here $\T$ acts on $S(\C)_+$ (freely off the base point) and $X$, but not
on $S^t$ or $H$.

The condition that $x \in E^2_{0,t} = H_t(X; \F_p)$ survives to the
$E^{2r}$-term, i.e., that all differentials $d^2(x), \dots, d^{2r-2}(x)$
vanish, is equivalent to $x$ being in the image from $H_t(F(S(\C^r)_+,
X)^{\T}; \F_p)$ under the map induced by restriction along $S(\C)_+
\subset S(\C^r)_+$.  This is in turn equivalent to the existence of a
$\T$-equivariant extension
$$
x' \colon S(\C^r)_+ \wedge S^t \to H \wedge X
$$
of $x$ along $S(\C)_+ \subset S(\C^r)_+$, in view of the natural
naively equivariant weak equivalence
$$
\nu \colon H \wedge F(S(\C^r)_+, X) \to F(S(\C^r)_+, H \wedge X) \,.
$$
(To establish this equivalence, use that $S(\C^r)_+$ is a finite $\T$-CW
complex.  We are considering maps from free $\T$-CW complexes into these
spectra, so only the naive notion of a $\T$-equivariant equivalence
is required.)

Suppose that $x \in E^{2r}_{0,t}$ has survived to the $E^{2r}$-term,
so that such a $\T$-equivariant extension $x'$ exists.  Then by the
construction of the homological spectral sequence, the differential
$$
d^{2r}(x) \in E^{2r}_{-2r, t+2r-1}
$$
is the obstruction to extending $x'$ further along $S(\C^r)_+ \subset
S(\C^{r+1})_+$ to an equivariant map
$$
x'' \colon S(\C^{r+1})_+ \wedge S^t \to H \wedge X \,.
$$
We put the obvious right adjoints of these maps together in a diagram.
$$
\xymatrix{
& S(\C)_+
\ar[dr]^-x \ar[d] \\
(\T \times \partial D^{2r})_+ \ar[r]^-{\alpha_+} \ar[d] &
S(\C^r)_+ \ar[r]^-{x'} \ar[d] & F(S^t, H \wedge X) \\
(\T \times D^{2r})_+ \ar[r] & S(\C^{r+1})_+
\ar@{-->}[ur]_-{x''}
}
$$
Recalling from~(\ref{e2.1}) that $S(\C^{r+1})_+$ is obtained from $S(\C^r)_+$
by adjoining a free $\T$-cell along the action map $\alpha \colon \T \times
\partial D^{2r} \to S(\C^r)$, the obstruction to such an extension
is precisely the obstruction to extending the equivariant map $x'
\circ \alpha_+$ from $(\T \times \partial D^{2r})_+$ over $(\T \times
D^{2r})_+$.  By adjunction, this equals the obstruction to extending
the underlying non-equivariant map $x' \colon \partial D^{2r}_+ = S(\C^r)_+
\to F(S^t, H \wedge X)$ over $D^{2r}_+$.  In terms of the preferred
stable splitting
$$
\partial D^{2r}_+ \simeq S^{2r-1} \vee D^{2r}_+
$$
(induced by the pinch map from $\partial D^{2r}_+$ to $S^{2r-1}$
and the inclusion to $D^{2r}_+$) this equals the restriction of the
non-equivariant map $x'$ to the stable summand $S^{2r-1}$.  Its left
adjoint again is a map
$$
\delta x \colon S^{2r-1} \wedge S^t \to H \wedge X \,,
$$
which represents $d^{2r}(x)$.  We summarize:

\begin{lemma}
\label{l3.1}
Let $x \in E^{2r}_{0,t} \subset H_t(X; \F_p)$ be represented by a
$\T$-equivariant map $x \colon S(\C)_+ \wedge S^t \to H \wedge X$ that extends
to a $\T$-equivariant map $x' \colon S(\C^r)_+ \wedge S^t \to H \wedge X$.
Then $d^{2r}(x) = y^r \cdot \delta x$, where $\delta x \in H_{t+2r-1}(X;
\F_p)$ is represented by $x'$ considered as a non-equivariant map,
restricted to the stable summand $S^{2r-1} \wedge S^t$ of $S(\C^r)_+
\wedge S^t$.
\end{lemma}

The extended map $x'$ represents a class in the homology of $F(S(\C^r)_+,
X)^{\T}$, and considering $x'$ as a non-equivariant map amounts to
following the map
$$
\phi \colon F(S(\C^r)_+, X)^{\T} \to F(S^{2r-1}_+, X)
$$
that forgets the $\T$-equivariance.  There is a canonical map
$$
\nu \colon X \wedge DS^{2r-1}_+ \to F(S^{2r-1}_+, X)
$$
where $DS^{2r-1}_+ = F(S^{2r-1}_+, S)$ is the functional dual of
$S^{2r-1}_+$, which is a weak equivalence since $S^{2r-1}_+$ is a finite
CW complex.  Hence there is a natural isomorphism
$$
\nu \colon H_*(X; \F_p) \otimes H^{-*}(S^{2r-1}; \F_p)
	\xrightarrow{\cong} H_*(F(S^{2r-1}_+, X); \F_p) \,,
$$
where we have identified the spectrum homology $H_*(DS^{2r-1}_+; \F_p)$
with the space-level cohomology $H^{-*}(S^{2r-1}; \F_p)$.  We write
$H^*(S^{2r-1}; \F_p) = E(\iota_{2r-1})$, where $\iota_{2r-1}$ is the
canonical generator in degree $(2r-1)$ and $E(-)$ denotes the exterior
algebra.

\begin{prop}
\label{p3.2}
The composite map
\begin{align*}
H_*(F(S(\C^r)_+, X)^{\T}; \F_p) & \xrightarrow{\phi_*} 
H_*(F(S^{2r-1}_+, X); \F_p) \\
&\xleftarrow[\cong]{\nu_*} H_*(X; \F_p) \otimes H^{-*}(S^{2r-1}; \F_p)
\end{align*}
takes any class $x'$ that maps to $x \in E^{2r}_{0,t} \subset H_t(X;
\F_p)$ by the restriction map
$$
H_*(F(S(\C^r)_+, X)^{\T}; \F_p) \to H_*(F(S(\C)_+, X)^{\T}; \F_p) =
H_*(X; \F_p)
$$
to the sum
$$
(\nu^{-1}_* \phi_*)(x') = x \otimes 1 + \delta x \otimes \iota_{2r-1}
\,,
$$
where $d^{2r}(x) = y^r \cdot \delta x$ in $E^{2r}_{-2r, t+2r-1}$.
Suppressing the power of $y$ we may somewhat imprecisely write this
formula as
$$
\phi_*(x) = x \otimes 1 + d^{2r}(x) \otimes \iota_{2r-1} \,.
$$
\end{prop}

The case $r=1$ of this result says that $d^2(x) = y \cdot \sigma x$,
and follows e.g.~from \cite[3.3]{Ro98}.
 
\begin{proof}
This is really a corollary to Lemma~\ref{l3.1}, but for the observation that
the restriction of the non-equivariant $x'$ to the subspace $S^t \subset
S(\C^r)_+ \wedge S^t$ equals the restriction of the non-equivariant
$x$ to the same subspace $S^t \subset S(\C)_+ \wedge S^t$, which in
turn corresponds to $x \in E^{2r}_{0,t}$ under the identification
$H_*(F(S(\C)_+, X)^{\T}; \F_p) = H_*(X; \F_p)$.  There are no signs
in these formulas, because the canonical map $\nu$ is derived from
the non-symmetric part of the closed monoidal structure on the stable
homotopy category.  (See e.g.~\cite[III.1]{LMS86}.)
\end{proof}

\begin{remark}
\label{r3.3}
Lemma~\ref{l3.1} says that the differential in the homotopy fixed point spectral
sequence is essentially the $\T$-equivariant root invariant for $H \wedge
X$.  A corresponding description of the Mahowald $C_2$-equivariant root
invariant for~$S$ can be found in \cite[2.5]{BG95}: Let $S^{n+k\alpha}$
denote the $C_2$-equivariant sphere that is the one point compactification
of $\R^n \oplus \R^k(-1)$, where $C_2$ acts trivially on $\R^n$ and by
negation on $\R^k(-1)$.  Given a non-equivariant (stable) map $x : S^n \to
S^0$, let $x' : S^{n+k\alpha} \to S^0$ be a $C_2$-equivariant extension
of $x$ with $k$ maximal.  Then the $C_2$-equivariant root invariant of $x$
contains the non-equivariant map $x' : S^{n+k} \to S^0$ underlying~$x'$.
\end{remark}

\section{ Commutative $S$-algebras }

Now suppose that $X = R$ is a naively $\T$-equivariant commutative
$S$-algebra, i.e., a commutative $S$-algebra with a continuous
point-set level action by the circle group~$\T$, through commutative
$S$-algebra maps.  We shall be concerned with the homotopy fixed points
of $R$, rather than its genuine fixed points, so only this weak notion
of an equivariant spectrum will be needed.  Cf.~\cite[Ch.~1]{GM95}.
Our principal example is $R = \THH(B)$, the topological Hochschild homology
spectrum of a non-equivariant commutative $S$-algebra $B$.  The cyclic
structure on topological Hochschild homology \cite[Ch.~IX]{EKMM97} then
provides the relevant $\T$-action.  In terms of the tensored structure
on commutative $S$-algebras over topological spaces, one model for the
topological Hochschild homology spectrum is $\THH(B) = B \otimes \T$,
and then $\T$ acts naturally through commutative $S$-algebra maps by
group multiplication in the $\T$-factor.

In this situation the homotopy fixed point spectrum $R^{h\T} = F(E\T_+,
R)^{\T}$ is also a commutative $S$-algebra.  Writing $\mu \colon R \wedge
R \to R$ for the $\T$-equivariant multiplication map of $R$, the
corresponding multiplication map for $R^{h\T}$ is given by the composite
\begin{align*}
F(E\T_+, R)^{\T} \wedge F(E\T_+, R)^{\T} &\xrightarrow{\wedge}
F((E\T \times E\T)_+, R \wedge R)^{\T} \\
&\xrightarrow{\mu_\# \Delta_+^\#} F(E\T_+, R)^{\T} \,.
\end{align*}
Here $\wedge$ smashes together two $\T$-equivariant maps $\Sigma^\infty
E\T_+ \to R$, and considers the resulting $(\T \times \T)$-equivariant map
as being $\T$-equivariant via the diagonal action.  The map $\mu_\#$
composes on the left by $\mu \colon R \wedge R \to R$, while the map
$\Delta_+^\#$ composes on the right by the diagonal map $\Delta_+ \colon
E\T_+ \to (E\T \times E\T)_+$.  Since $\mu$ is commutative and $\Delta_+$
is cocommutative, the resulting multiplication on $R^{h\T}$ is also
strictly commutative.

Writing $\eta \colon S \to R$ for the $\T$-equivariant unit map of
$R$, the corresponding unit map for $R^{h\T}$ is the composite
$$
S \to F(E\T_+, S)^{\T} \xrightarrow{\eta_\#} F(E\T_+, R)^{\T} \,.
$$
To define the first map, we must use that $\T$ acts trivially on~$S$.

The same constructions can be applied for the $\T$-CW skeleta of $E\T$,
so each $F(S(\C^n)_+, R)^{\T}$ is also a commutative $S$-algebra.
In particular, the continuous homology $H^c_*(R^{h\T}; \F_p)$ is a
limit of graded commutative algebras, and is therefore itself a graded
commutative algebra.

\begin{prop}
\label{p4.1}
Let $R$ be a $\T$-equivariant commutative $S$-algebra.  Then
the homological homotopy fixed point spectral sequence
$$
E^2_{**}(R) = P(y) \otimes H_*(R; \F_p)
\Longrightarrow H^c_*(R^{h\T}; \F_p)
$$
is an $A_*$-comodule algebra spectral sequence, where $P(y)$ is the
polynomial algebra on $y$ in bidegree $(-2,0)$, and $H_*(R; \F_p)$
has the Pontryagin product.
\end{prop}

\proof
The algebra product in $E^r(R)$ is derived from the $\T$-equivariant
commutative $S$-algebra product in $R$ and a $\T$-equivariant cellular
approximation~$d$ to the diagonal map $\Delta \colon E\T \to E\T \times E\T$.
The deduction is in principle standard, but due to our homological
indexing and perhaps unusual choice of exact couple defining the
spectral sequence, it is not so easy to find an applicable reference.
(The closest may be a combination of Remark~\ref{r2.7}(c) and 
\cite[4.3.5]{HM03}, adapted from the $G$-Tate construction for a finite group $G$
to the $\T$-homotopy fixed point spectrum.)  We therefore provide the
following outline.

The product $(\T \times \T)$-CW structure on $E\T \times E\T$ can be
refined to a $\T$-CW structure, by starting with a $\T$-CW structure
on $\T \times \T$ with the diagonal $\T$-action.  Fix a choice of a
$\T$-equivariant cellular map $d \colon E\T \to E\T \times E\T$ that is
equivariantly homotopic to $\Delta$.

To shorten the notation, we shall write
$$
E\T^m_a = E\T^{(m)}/E\T^{(a-1)}
$$
within this proof.  Then for all integers $s$, $s'$ and $r\ge1$, $d$
induces a map of subquotients
$$
d \colon {E\T^{-s-s'+r-1}_{-s-s'}}
\to {E\T^{-s+r-1}_{-s}} \wedge {E\T^{-s'+r-1}_{-s'}} \,,
$$
since each $\T$-$m$-cell in the source maps into the $(\T \times
\T)$-$m$-skeleton in the product structure on the target, which is the
union of products of $\T$-$k$-cells and $\T$-$k'$-cells with $k \ge -s$,
$k' \ge -s'$ and $k+k' \le m$.  For $s = -2n$ and $s' = -2n'$ non-positive
and even, and $r=1$, this is a map
\begin{equation}
d \colon \T_+ \wedge S^{2n+2n'} \to \T_+ \wedge S^{2n} \wedge \T_+ \wedge S^{2n'}
\label{e4.2}
\end{equation}
homotopic to the $(2n+2n')$-th suspension of the diagonal $\T_+ \to
(\T \times \T)_+$.

For integers $a \le b$, let
$$
H(a,b) = H_*(F(E\T^{b-1}_a, R)^{\T}; \F_p)
$$
be the finite terms of a Cartan--Eilenberg system \cite[XV.7]{CE56},
with $E^1_s = H(-s,-s+1)$.  Applying $F(-, R)^{\T}$ and homology to $d$,
we obtain homomorphisms
$$
\phi_r \colon H(-s,-s+r) \otimes H(-s',-s'+r) \to H(-s-s',-s-s'+r)
$$
for all integers $s$, $s'$ and $r\ge1$.  For $s = -2n$, $s' = -2n'$
non-positive and $r=1$, this is the homomorphism
$$
\phi_1 \colon E^1_{-2n} \otimes E^1_{-2n'} \to E^1_{-2n-2n'}
$$
that takes $y^n \cdot x \otimes y^{n'} \cdot x'$ to $y^{n+n'} \cdot
\mu_*(x \otimes x')$, under the identification $E^1_{-2n,t} = y^n
\cdot H_t(R; \F_p)$ from Section~2.  So $E^1 = E^2$ equipped with
the product $\phi_1$ is isomorphic to the tensor product of the polynomial
algebra $P(y)$ and $H_*(R; \F_p)$ with the Pontryagin product $\mu_*$.

To verify that the spectral sequence differentials $d^r$ are derivations,
so that each term $E^{r+1}$ inductively inherits an algebra structure
from $E^r$, it suffices to check that we have a multiplicative
Cartan--Eilenberg system, i.e., that the pairings $\phi_r$ satisfy
the relation
\begin{equation}
\delta(\phi_r(z \otimes z')) = \phi_1(\delta(z) \otimes \eta(z'))
+ (-1)^{|z|} \phi_1(\eta(z) \otimes \delta(z'))
\label{e4.3}
\end{equation}
in $E^1_{s+s'-r}$, for $z \in H(-s,-s+r)$ of degree~$|z|$ and $z'
\in H(-s',-s'+r)$.  Here $\delta \colon H(-s,-s+r) \to H(-s+r,-s+r+1) =
E^1_{s-r}$ is the degree $(-1)$ homomorphism induced by the stable
connecting map
$$
\delta \colon \Sigma^{-1} {E\T^{-s+r}_{-s+r}}
\to {E\T^{-s+r-1}_{-s}}
$$
of the obvious triple, and $\eta \colon H(-s,-s+r) \to H(-s,-s+1) =
E^1_s$ is the homomorphism induced by the inclusion
$$
\eta \colon {E\T^{-s}_{-s}}
\to {E\T^{-s+r-1}_{-s}} \,.
$$
We use similar notations with $s'$ and $s+s'$ in place of $s$.
The sufficiency of this condition can be read directly off from the
definition of the differentials in the spectral sequence associated to
an unrolled exact couple, or obtained from \cite{Ma54} or 
\cite[Ex.~2.2.1]{Mc01}.

In geometric terms, (\ref{e4.3}) asks that the composite map
$$
\Sigma^{-1} {E\T^{-s-s'+r}_{-s-s'+r}}
\xrightarrow{\delta} 
{E\T^{-s-s'+r-1}_{-s-s'}}
\xrightarrow{d} {E\T^{-s+r-1}_{-s}}
\wedge {E\T^{-s'+r-1}_{-s'}}
$$
is homotopic, as a map of $\T$-equivariant spectra, to the sum of the
two composite maps
$$
\Sigma^{-1} {E\T^{-s-s'+r}_{-s-s'+r}}
\xrightarrow{d} \Sigma^{-1} {E\T^{-s+r}_{-s+r}}
\wedge {E\T^{-s'}_{-s'}}
\xrightarrow{\delta \wedge \eta}
{E\T^{-s+r-1}_{-s}}
\wedge {E\T^{-s'+r-1}_{-s'}}
$$
and
$$
\Sigma^{-1} {E\T^{-s-s'+r}_{-s-s'+r}}
\xrightarrow{(-1)^s d} {E\T^{-s}_{-s}}
\wedge \Sigma^{-1} {E\T^{-s'+r}_{-s'+r}}
\xrightarrow{\eta \wedge \delta}
{E\T^{-s+r-1}_{-s}}
\wedge {E\T^{-s'+r-1}_{-s'}}
\,.
$$

There is only something to check when $s$, $s'$ and $r$ are all even,
with $s$ and $s'$ non-positive.  The common source of the $\T$-equivariant
stable maps to be compared is then $\Sigma^{-1} \T_+ \wedge S^{-s-s'+r}$,
so by an adjunction we may as well compare non-equivariant maps from
$S^{-s-s'+r-1}$ to
$$
E\T^{-s+r-1}_{-s} \wedge E\T^{-s'+r-1}_{-s'}
\simeq (S^{-s+r-1} \vee S^{-s}) \wedge (S^{-s'+r-1} \vee S^{-s'}) \,.
$$
The projection to the stable summand $\Sigma E\T^{(-s-1)} \wedge \Sigma
E\T^{(-s'-1)} = S^{-s} \wedge S^{-s'}$ in the target is trivial for each
of the three maps, since in each case two subsequent maps in a cofiber
sequence occur as a factor of the composite map.  The projections to the
summands $E\T^{(-s+r-1)} \wedge \Sigma E\T^{(-s'-1)} = S^{-s+r-1} \wedge
S^{-s'}$ and $\Sigma E\T^{(-s-1)} \wedge E\T^{(-s'+r-1)} = S^{-s} \wedge
S^{-s'+r-1}$ agree as required, by the same kind of homotopy as in~(\ref{e4.2}).

Finally, in the presence of strong convergence due to the vanishing
of Boardman's obstruction group $RE^\infty_{**}$, we claim that the
isomorphisms $F_s/F_{s-1} \cong E^\infty_s$ take the associated graded
algebra structure derived from the product on $H^c_*(R^{h\T}; \F_p)$
to the algebra structure on $E^\infty$-term.  These isomorphisms are
obtained by descent to subquotients, as in~(\ref{e2.4}), from the homomorphisms
$$
\lim_m
H_*(F(E\T^m_{-s}, R)^{\T}; \F_p) = H(-s,\infty) \to H(-s,-s+1) = E^1_s \,,
$$
so it suffices to verify that the following diagram commutes:
$$
\xymatrix{
H(-s,\infty) \otimes H(-s',\infty) \ar[r] \ar[d] 
	& H(-s-s',\infty) \ar[d] \\
E^1_s \otimes E^1_{s'} \ar[r]^-{\phi_1} & E^1_{s+s'} 
}
$$
This follows immediately from the commutativity of the following diagram,
%% of $\T$-CW complexes, 
where $d$ is the diagonal approximation and the
vertical maps are inclusions:
$$
\xymatrix{
E\T/E\T^{(-s-1)} \wedge E\T/E\T^{(-s'-1)} &
E\T/E\T^{(-s-s'-1)} \ar[l]_-d \\
E\T^{(-s)}/E\T^{(-s-1)} \wedge E\T^{(-s')}/E\T^{(-s'-1)} \ar[u] &
E\T^{(-s-s')}/E\T^{(-s-s'-1)} \ar[u] \ar[l]_-d \\
}\eqno{\smash{\raise -45pt\hbox{\qed}}}
$$

Commutative $S$-algebras are $E_\infty$ ring spectra 
(\cite[II.4]{EKMM97}), and are in particular also $H_\infty$ ring spectra 
(\cite[III.5]{EKMM97} and \cite[VII.2]{LMS86}).  Hence there are Dyer--Lashof operations
$Q^i$ acting on their mod~$p$ homology algebras \cite[III.1]{BMMS86}.
Recall that $Q^i$ is a natural transformation
$$
Q^i \colon H_t(R; \F_p) \to H_{t+iq}(R; \F_p)
$$
for all integers $t$, where $q = 2p-2$.  We also include their composites
$\beta Q^i$ with the homology Bockstein operation $\beta \colon H_t(R;
\F_p) \to H_{t-1}(R; \F_p)$ as generators of the Dyer--Lashof algebra.
For $p=2$ the standard notation is to write $Q^{2i}$ and $Q^{2i-1}$
for the operations that would otherwise be called $Q^i$ and $\beta Q^i$,
respectively.

As noted before Proposition~\ref{p4.1}, the commutative $S$-algebra structure
on $R^{h\T} = F(E\T_+, R)^{\T}$ restricts to one on each $F(S(\C^n)_+,
R)^{\T}$, so each algebra homomorphism
$$
H^c_*(R^{h\T}; \F_p) \to H_*(F(S(\C^n)_+, R)^{\T}; \F_p)
$$
commutes with the Dyer--Lashof operations in the source and in the target.
In particular, for $n=1$ the identification $F(S(\C)_+, R)^{\T} \cong R$
lets us recognize the Dyer--Lashof operations in the target as those in
$H_*(R; \F_p)$.  This action by Dyer--Lashof operations on the vertical
axis $E^2_{0,*}$ can be algebraically extended to an action on the full
$E^2$-term, by the formula $\beta^\epsilon Q^i(y^n \cdot x) = y^n \cdot
\beta^\epsilon Q^i(x)$.  Our next result shows that this action extends
further to all terms of the homological homotopy fixed point spectral
sequence.

\begin{prop}
\label{p4.4}
$R$ be a $\T$-equivariant commutative $S$-algebra and let $E^r(R)$ be
its homological homotopy fixed point spectral sequence.  Then for each
element $x \in E^{2r}_{0,t} \subset H_t(R; \F_p)$ we have the relation
$$
d^{2r}(\beta^\epsilon Q^i(x)) = \beta^\epsilon Q^i(d^{2r}(x)) \,,
$$
for every integer $i$ and $\epsilon \in \{0,1\}$.  Here the right
hand side should be interpreted as follows:  If $d^{2r}(x) = y^r \cdot
\delta x$ with $\delta x \in H_{t+2r-1}(R; \F_p)$ then $\beta^\epsilon
Q^i(d^{2r}(x)) = y^r \cdot \beta^\epsilon Q^i(\delta x)$.
\end{prop}

The case $r=1$ appears as \cite[5.9]{AnR}, with a proof that we generalize
as follows.

\begin{proof}
Let $x \in H_t(R; \F_p)$ and suppose that $x$ survives to the
$E^{2r}$-term.  Then there exists an extension $x' \in H_t(F(S(\C^r)_+,
R)^{\T}; \F_p)$ of $x$ over the restriction map
$$
F(S(\C^r)_+, R)^{\T} \to F(S(\C)_+, R)^{\T} \cong R
$$
of commutative $S$-algebras, and $z' = \beta^\epsilon Q^i(x')$ is
an extension of $z = \beta^\epsilon Q^i(x)$ over the same map, by
naturality of the Dyer--Lashof operations.  The maps $\phi$ and $\nu$
from Proposition~\ref{p3.2} are both maps of commutative $S$-algebras, and
therefore induce algebra homomorphisms $\phi_*$ and $\nu_*$ that commute
with the Dyer--Lashof operations.  Thus
\begin{equation}
(\nu_*^{-1} \phi_*)(\beta^\epsilon Q^i(x')) = \beta^\epsilon
Q^i(x) \otimes 1 + \delta z \otimes \iota_{2r-1}
\label{e4.5}
\end{equation}
where $d^{2r}(\beta^\epsilon Q^i(x)) = y^r \cdot \delta z$, is equal to
$$
\beta^\epsilon Q^i((\nu_*^{-1} \phi_*)(x')) =
\beta^\epsilon Q^i(x \otimes 1 + \delta x \otimes \iota_{2r-1})
$$
where $d^{2r}(x) = y^r \cdot \delta x$.  Now the Dyer--Lashof operations
on the homology of the smash product $R \wedge DS^{2r-1}_+$ are given by
a Cartan formula, and on the tensor factor $H_*(DS^{2r-1}_+; \F_p) \cong
H^{-*}(S^{2r-1}; \F_p)$ the operation $\beta^\epsilon Q^i$ corresponds to
the Steenrod operation $\beta^\epsilon P^{-i}$, by \cite[III.1.2]{BMMS86}.
But the latter operations all act trivially on $H^*(S^{2r-1}; \F_p)$,
except for $P^0 = 1$, so the Cartan formula gives
$$
\beta^\epsilon Q^i(x \otimes 1 + \delta x \otimes \iota_{2r-1})
= \beta^\epsilon Q^i(x) \otimes 1 +
	\beta^\epsilon Q^i(\delta x) \otimes \iota_{2r-1} \,.
$$
Identifying this with~(\ref{e4.5}) and comparing the coefficients
of $\iota_{2r-1}$ we obtain the identity
$$
\delta z = \beta^\epsilon Q^i(\delta x) \,,
$$
as claimed.
\end{proof}

The homological homotopy fixed point spectral sequence for $S$ itself
is particularly simple.  For $H_*(S; \F_p) = \F_p$, so the spectral
sequence collapses to
$$
E^2_{**}(S) = P(y) \,,
$$
concentrated on the horizontal axis.  Hence each power of $y$ is
an infinite cycle, i.e., $d^{2r}(y^n) = 0$ for all $r\ge1$ and $n$.
The spectral sequence $E^*(X)$ for a general $\T$-equivariant $S$-module
$X$ is a module over the one for $S$, so the Leibniz formula for the
module pairing immediately yields part~(a) of the following result.

\begin{lemma}
\label{l4.6}
Let $X$ be any $\T$-equivariant $S$-module.

{\rm(a)}\qua
The differentials in the homological homotopy fixed point spectral
sequence converging to $H^c_*(X^{h\T}; \F_p)$ satisfy the relation
$$
d^{2r}(y^n \cdot x) = y^n \cdot d^{2r}(x)
$$
for all $x \in E^{2r}_{0,*} \subset H_*(X; \F_p)$, $2r\ge2$
and~$n\ge0$.

{\rm(b)}\qua
Each class in $E^{2r}_{-2n,t}$ has the form $y^n \cdot x$ for a
class $x \in E^{2r}_{0,t} \subset H_t(X; \F_p)$.  Hence the spectral
sequence is completely determined by the differentials that originate
on the vertical axis.

{\rm(c)}\qua
The $y$-torsion in $E^{2r}_{**}$ has height strictly less than~$r$,
and is concentrated in filtrations $-2r+4 \le s \le 0$.
\end{lemma}

\begin{proof}
It remains to prove parts~(b) and~(c), which we do by induction on $r$.
The claims for $r=1$ are clear from Propositions~\ref{p2.6} and~\ref{p4.1}.

If $z = y^n \cdot x \in E^{2r}_{-2n,t}$ survives to the
$E^{2r+2}$-term, then $0 = d^{2r}(y^n \cdot x) = y^n \cdot d^{2r}(x)$
by part~(a).  By induction there is no $y$-torsion in or below the
filtration of $d^{2r}(x)$, so $d^{2r}(x) = 0$, $x$ survives to the
$E^{2r+2}$-term, and we still have $z = y^n \cdot x$.

If now $y \cdot z = 0$ in the $E^{2r+2}$-term, then $y \cdot z$ in the
$E^{2r}$-term must be a boundary of the form $d^{2r}(y^{n+1-r} \cdot w)$,
with $w$ on the vertical axis.  If $n \ge r$ then it follows that $z =
d^{2r}(y^{n-r} \cdot w)$ in the $E^{2r}$-term, since by induction there
is no $y$-torsion in filtration $-2n$ of the $E^{2r}$-term.  Thus the
only $y$-torsion in the $E^{2r+2}$-term lies in filtrations $-2r+2 \le
s \le 0$.
\end{proof}

In Section~7 we shall remark on an analogous homological Tate spectral
sequence, where $P(y)$ is replaced by $P(y, y^{-1})$ and the issue of
$y$-torsion classes becomes void.

\section{ Infinite cycles }

The Dyer--Lashof operations satisfy instability conditions 
\cite[III.1.1]{BMMS86} that are in a sense dual to those of the Steenrod operations.
For a class $x \in H_t(R; \F_p)$ the lowest nontrivial operation is
$Q^t(x) = x^2$ when $p=2$, $Q^m(x) = x^p$ when $p$ is odd and $t =
2m$ is even, and $\beta Q^m(x)$ when $p$ is odd and $t = 2m-1$ is odd.
Similarly, the lowest nontrivial operation on $\delta x \in H_{t+2r-1}(R;
\F_p)$ with $d^{2r}(x) = y^r \cdot \delta x$ is $Q^{t+2r-1}(\delta x) =
(\delta x)^2$ when $p=2$, $\beta Q^{m+r}(\delta x)$ when $p$ is odd and
$t = 2m$ is even, and $Q^{m+r-1}(\delta x) = (\delta x)^p$ when $p$
is odd and $t = 2m-1$ is odd.  Thus there is in each case a sequence
of $(2r-1)$ Dyer--Lashof operations $\beta^\epsilon Q^i$ whose
action on $x$ can be nontrivial, but whose action on $\delta x$ must
be trivial.  By Proposition~\ref{p4.4}, this sequence of operations on $x$
will survive past the $E^{2r}$-term, at least to the $E^{2r+2}$-term.
It is the main point of the present article to show that these classes,
and one more ``companion class'', then in fact go on indefinitely to
survive to the $E^\infty$-term, i.e., are infinite cycles!

\begin{thm}
\label{t5.1}
Let $R$ be a $\T$-equivariant commutative $S$-algebra, suppose that
$x \in H_t(R; \F_p)$ survives to the $E^{2r}$-term $E^{2r}_{0,t}
\subset H_t(R; \F_p)$ of the homological homotopy fixed point spectral
sequence for $R$, and write $d^{2r}(x) = y^r \cdot \delta x$.

{\rm(a)}\qua
For $p=2$, the $2r$ classes
$$
x^2 = Q^t(x),\  Q^{t+1}(x),\  \dots \ ,\  Q^{t+2r-2}(x)
\quad\text{and}\quad Q^{t+2r-1}(x) + x \delta x
$$
all survive to the $E^\infty$-term, i.e., are infinite cycles.

{\rm(b)}\qua
For $p$ odd and $t = 2m$ even, the $2r$ classes
$$
x^p = Q^m(x),\  \beta Q^{m+1}(x),\  \dots \ ,\  Q^{m+r-1}(x)
\quad\text{and}\quad x^{p-1} \delta x
$$
all survive to the $E^\infty$-term, i.e., are infinite cycles.

{\rm(c)}\qua
For $p$ odd and $t = 2m-1$ odd, the $2r$ classes
$$
\beta Q^m(x),\  Q^m(x),\  \dots \ ,\  \beta Q^{m+r-1}(x)
\quad\text{and}\quad Q^{m+r-1}(x) - x(\delta x)^{p-1}
$$
all survive to the $E^\infty$-term, i.e., are infinite cycles.
\end{thm}

\begin{proof}
The argument proceeds by considering a universal example.  Recall that
a class $x \in E^{2r}_{0,t}$ is represented by a $\T$-equivariant
map $x \colon S(\C)_+ \wedge S^t \to H \wedge R$ that admits an equivariant
extension $x' \colon S(\C^r)_+ \wedge S^t \to H \wedge R$.  Let
$$
X = D_p(S(\C^r)_+ \wedge S^t) = E\Sigma_p \ltimes_{\Sigma_p} (S(\C^r)_+
\wedge S^t)^{\wedge p}
$$
be the $p$-th extended power of the spectrum $S(\C^r)_+ \wedge S^t$.

Somewhat abusively, we write $H_*(S(\C^r)_+ \wedge S^t; \F_p)
= \F_p\{x, \delta x\}$ with $|x| = t$ and $|\delta x| = t+2r-1$.
Then the homology of the $p$-th extended power spectrum is
$$
H_*(X; \F_2) = \F_2\{x \delta x,\,\,  Q^i(x) \mid i\ge t,\,\, Q^i(\delta x)
\mid i\ge t+2r-1\}
$$
for $p=2$,
$$
H_*(X; \F_p) = \F_p\{x^{p-1} \delta x,\,\, \beta^\epsilon Q^i(x) \mid i \ge
m+\epsilon,\,\, \beta^\epsilon Q^i(\delta x) \mid i \ge m+r\}
$$
for $p$ odd and $t = 2m$ even, and
$$
H_*(X; \F_p) = \F_p\{x (\delta x)^{p-1},\,\, \beta^\epsilon Q^i(x) \mid i
\ge m,\,\, \beta^\epsilon Q^i(\delta x) \mid i \ge m+r-1+\epsilon\}
$$
for $p$ odd and $t = 2m-1$ odd.  Throughout $i$ is an integer and
$\epsilon \in \{0,1\}$.

The equivariant extension $x'$ induces an equivariant map
$$
D_p(x') \colon X = D_p(S(\C^r)_+ \wedge S^t) \to D_p(H \wedge R) \,.
$$
The $\T$-equivariant commutative $S$-algebra structures on $H$ and $R$
combine to form one on $H \wedge R$, and, as noted in Section~4, this
gives $H \wedge R$ an $H_\infty$-ring structure 
(\cite[Thm.~0.1, II.4 and III.5]{EKMM97}
and \cite[VII.2]{LMS86}).  The associated
$H_\infty$ structure includes, in particular, a $\T$-equivariant
structure map
$$
\xi_p \colon D_p(H \wedge R) \to H \wedge R
$$
that extends the $p$-fold multiplication map on $H \wedge R$.
Taken together, these produce an equivariant map
$$
H \wedge D_p(S(\C^r)_+ \wedge S^t)
\xrightarrow{1\wedge D_p(x')} H \wedge D_p(H \wedge R) \xrightarrow{1\wedge \xi_p}
H \wedge H \wedge R \xrightarrow{\mu\wedge 1} H \wedge R \,,
$$
where $\mu$ is the multiplication on $H$.  Applying homotopy we have
a homomorphism
\begin{equation}
H_*(X; \F_p) = H_*(D_p(S(\C^r)_+ \wedge S^t); \F_p)
\to H_*(R; \F_p)
\label{e5.2}
\end{equation}
which, by definition, takes the classes generating $H_*(X; \F_p)$ to the
classes with the same names in $H_*(R; \F_p)$.
Now $X = D_p(S(\C^r)_+ \wedge S^t)$ is a $\T$-equivariant retract of
the free commutative $S$-algebra
$$
P \simeq \bigvee_{j\ge0} D_j(S(\C^r)_+ \wedge S^t)
$$
on the space $S(\C^r)_+ \wedge S^t$, so the homological homotopy fixed
point spectral sequence for $X$ is a direct summand of the one for $P$.
Thus the formula from Proposition~\ref{p4.4} for the $d^{2r}$-differentials
in the spectral sequence for $P$ is also applicable in the spectral
sequence for $X$.

Now consider the homological homotopy fixed point spectral sequence for
$X = D_p(S(\C^r)_+ \wedge S^t)$, first for $p=2$ and then for $p$ odd.
We shall show in each case that the $2r$ classes in $E^{2r}_{0,*}
\subset H_*(X; \F_p)$, with names as listed in the statement of the
theorem, are infinite cycles.  By naturality of the homotopy fixed
point spectral sequence with respect to the map $H \wedge X \to H
\wedge R$ from~(\ref{e5.2}), it follows that the $2r$ target classes listed
in $E^{2r}_{0,*} \subset H_*(R; \F_p)$ are also infinite cycles.
This will complete the proof of the theorem.

{\rm(a)}\qua
Let $p=2$.  The homological homotopy fixed point spectral sequence for
$X$ has
$$
E^2_{**} = P(y) \otimes \F_2\{x \delta x,\,\, Q^i(x) \mid i\ge t,\,\,
Q^i(\delta x) \mid i \ge t+2r-1\}
$$
and nontrivial differentials $d^{2r}(x \delta x) = y^r \cdot (\delta
x)^2$ and
$$
d^{2r}(Q^i(x)) = y^r \cdot Q^i(\delta x)
$$
for all $i\ge t+2r-1$, together with their $y$-multiples.  

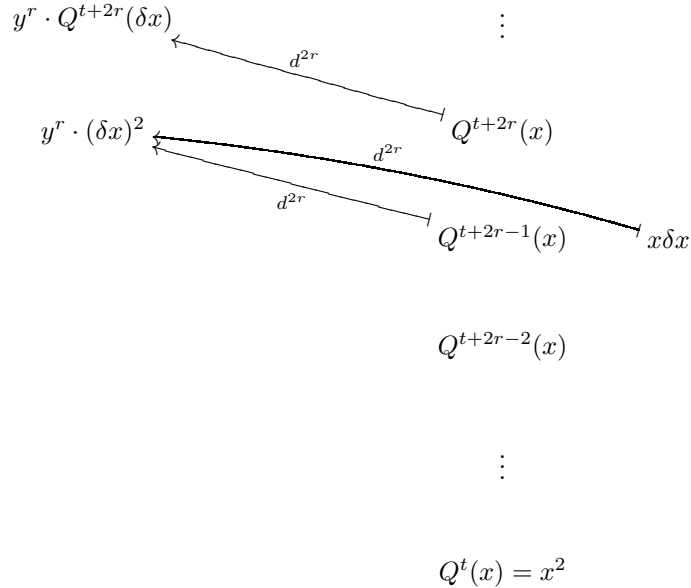
\begin{figure}[ht!]\small
$$
\xymatrix{
y^r \cdot Q^{t+2r}(\delta x) && \vdots \\
y^r \cdot (\delta x)^2 &&
	Q^{t+2r}(x) \ar@{|->}[ull]_{d^{2r}} \\
&& Q^{t+2r-1}(x) \ar@{|->}[ull]^{d^{2r}}
	& x \delta x \ar@{|->}@/_/[ulll]_{d^{2r}} \\
&& Q^{t+2r-2}(x) \\
&& \vdots \\
& \qquad\qquad & Q^t(x) = x^2 \\
}
$$
\caption{The case $p=2$}
\end{figure}

This leaves
$$
E^{2r+2}_{**} = P(y) \otimes \F_2\{Q^i(x) \mid t \le i < t+2r-1,\,\,
Q^{t+2r-1}(x) + x \delta x\}
$$
plus some $y$-torsion classes from $E^2_{**}$ in filtrations $-2r <
s \le 0$.  Hence there are no classes remaining in the entire quadrant
with filtration $s \le -2r$ and vertical degree $* > |x \delta x| =
2t+2r-1$.  All further differentials on the classes in $E^{2r+2}_{0,*}$
on the vertical axis land in this zero region, since already $E^2_{0,*}$
starts in degree $2t$ with the lowest class $Q^t(x) = x^2$.  Thus all
further differentials from the vertical axis are zero, and the spectral
sequence collapses at $E^{2r+2}_{**} = E^\infty_{**}$.

{\rm(b)}\qua
Let $p$ be odd and $t = 2m$ even.  The homological homotopy fixed
point spectral sequence for $X$ has
$$
E^2_{**} = P(y) \otimes \F_p\{x^{p-1} \delta x,\,\, \beta^\epsilon Q^i(x)
\mid i \ge m+\epsilon,\,\, \beta^\epsilon Q^i(\delta x) \mid i \ge m+r\}
$$
and nontrivial differentials
$$
d^{2r}(\beta^\epsilon Q^i(x)) = y^r \cdot \beta^\epsilon Q^i(\delta x)
$$
for $i \ge m+r$.  This leaves
$$
E^{2r+2}_{**} = P(y) \otimes \F_p\{x^{p-1}\delta x,\,\, Q^i(x)
\mid m+\epsilon \le i < m+r\}
$$
plus some $y$-torsion classes in filtrations $-2r < s \le 0$.  Hence there
are no classes left in the region where $s \le -2r$ and the vertical
degree is $* > |Q^{m+r-1}(x)|$.

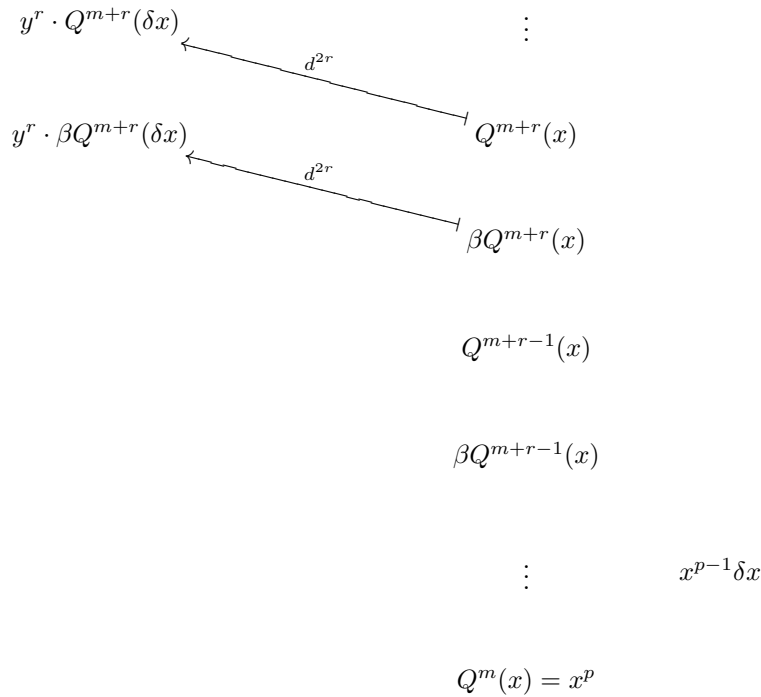
\begin{figure}[ht!]\small
$$
\xymatrix{
y^r \cdot Q^{m+r}(\delta x) && \vdots \\
y^r \cdot \beta Q^{m+r}(\delta x) &&
	Q^{m+r}(x) \ar@{|->}[ull]_{d^{2r}} \\
&& \beta Q^{m+r}(x) \ar@{|->}[ull]_{d^{2r}} \\
&& Q^{m+r-1}(x) \\
&& \beta Q^{m+r-1}(x) \\
&& \vdots & x^{p-1} \delta x \\
& \qquad\qquad & Q^m(x) = x^p \\
}
$$
\medskip
\caption{The case $p$ odd and $t=2m$ even}
\end{figure}

Now, $x$ was also a class in the $E^{2r-2}$-term, with $d^{2r-2}(x) = 0$,
so by induction over $r$ we may assume (by naturality from the case of
$(r-1)$) that the classes $\beta^\epsilon Q^i(x)$ with $m+\epsilon \le i
< m+(r-1)$ are infinite cycles.  This leaves the three classes $x^{p-1}
\delta x$, $\beta Q^{m+r-1}(x)$ and $Q^{m+r-1}(x)$ in $E^{2r+2}_{0,*}$
that are not $y$-torsion, and could therefore imaginably support
a differential after $d^{2r}$.  But the first two classes $\beta
Q^{m+r-1}(x)$ and $Q^{m+r-1}(x)$ are so close to the horizontal edge of
the vanishing region that all differentials after $d^{2r}$ must vanish
on these classes.

The third class $x^{p-1} \delta x$ has odd degree, so an even length
differential on it must land in an even degree.  The only even degree
classes in filtrations $s \le -2r$ are the $y$-multiples of $Q^i(x)$ for
$m \le i < m+r$, of which $Q^m(x) = x^p$ is in lower degree than that
of $x^{p-1} \delta x$.  The remaining possible target classes $Q^i(x)$
for $m < i < m+r$ all have nontrivial Bockstein images $\beta Q^i(x)$,
but $\beta(x^{p-1} \delta x) = 0$ in $H_*(X; \F_p)$.  Therefore, by
naturality of the differential with respect to the Bockstein operation,
all of these targets for a differential on $x^{p-1} \delta x$ are
excluded.  Thus also $x^{p-1} \delta x$ is an infinite cycle.

{\rm(c)}\qua
Let $p$ be odd and $t = 2m-1$ odd.  The homological homotopy fixed
point spectral sequence for $X$ has
$$
E^2_{**} = P(y) \otimes \F_p\{x (\delta x)^{p-1},\,\, \beta^\epsilon Q^i(x)
\mid i \ge m,\,\, \beta^\epsilon Q^i(\delta x) \mid i \ge m+r-1+\epsilon\}
$$
and nontrivial differentials $d^{2r}(x (\delta x)^{p-1}) = y^r \cdot
(\delta x)^p$ and
$$
d^{2r}(\beta^\epsilon Q^i(x)) = y^r \cdot \beta^\epsilon Q^i(\delta x)
$$
for $i \ge m+r-1+\epsilon$.  This leaves
$$
E^{2r+2}_{**} = P(y) \otimes \F_p\{\beta^\epsilon Q^i(x)
\mid m \le i < m+r-1+\epsilon,\,\, Q^{m+r-1}(x) - x(\delta x)^{p-1}\}
$$
plus $y$-torsion classes in filtrations $-2r < s \le 0$.  Hence there are no
classes left in the region where $s \le -2r$ and the vertical degree
is $* > |Q^{m+r-1}(x)|$.

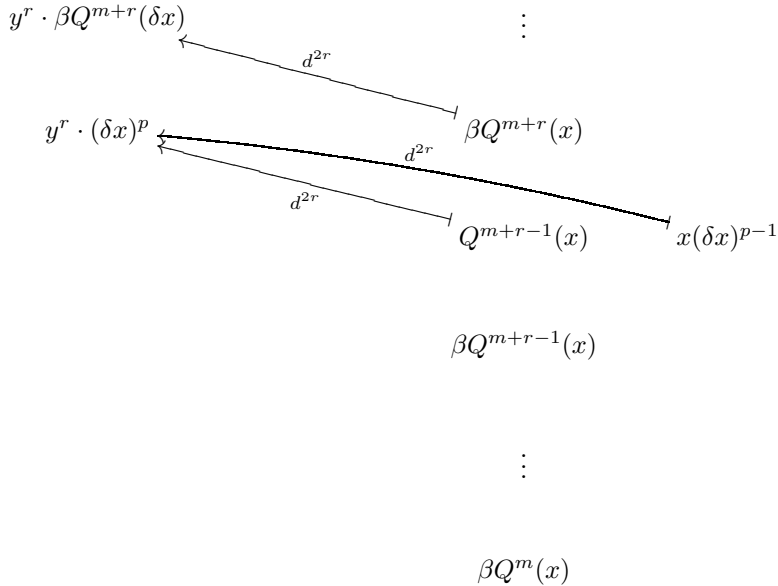
\begin{figure}[ht!]\small
$$
\xymatrix{
y^r \cdot \beta Q^{m+r}(\delta x) && \vdots \\
y^r \cdot (\delta x)^p &&
	\beta Q^{m+r}(x) \ar@{|->}[ull]_{d^{2r}} \\
&& Q^{m+r-1}(x) \ar@{|->}[ull]^{d^{2r}}
	& x(\delta x)^{p-1} \ar@{|->}@/_/[ulll]_{d^{2r}} \\
&& \beta Q^{m+r-1}(x) \\
&& \vdots \\
& \qquad\qquad & \beta Q^m(x) \\
}
$$
\caption{The case $p$ odd and $t=2m-1$ odd}
\end{figure}

Again considering $x$ as a class in $E^{2r-2}_{0,*}$ and using
induction on $r$ we may
assume that the classes $\beta^\epsilon Q^i(x)$ for $m \le i <
m+r-2+\epsilon$ and $Q^{m+r-2}(x) - x(\delta'x)^{p-1}$ are infinite
cycles.  Here $\delta'x$ is defined by $d^{2r-2}(x) = y^{r-1} \cdot
\delta'x$.  The fact that $d^{2r-2}(x) = 0$ gives
$\delta'x = 0$, so in fact all the classes $\beta^\epsilon Q^i(x)$
for $m \le i < m+r-1$ in $E^{2r+2}_{**}$ are infinite cycles.

This leaves only the two classes $\beta Q^{m+r-1}(x)$ and $Q^{m+r-1}(x)
- x(\delta x)^{p-1}$, but these are so close to the horizontal border
of the vanishing region that all differentials after $d^{2r}$ must be
zero on them.
\end{proof}

\section{ Examples }

Our Theorem~\ref{t5.1} has applications to the homological homotopy fixed point
spectral sequence for the commutative $S$-algebra $R = \THH(B)$ given by
the topological Hochschild homology of a commutative $S$-algebra $B$.
The $\T$-homotopy fixed point spectrum $\THH(B)^{h\T}$ is closely related
to the topological model $TF(B)$ for the negative cyclic homology of
$B$, which in turn is very close to the topological cyclic homology
$TC(B)$ \cite{BHM93} and algebraic $K$-theory $K(B)$ of $B$ \cite{Du97}.
These spectral sequences therefore have significant interest.

First consider the connective Johnson--Wilson spectrum $B = BP\langle
m{-}1\rangle$, for some prime $p$ and integer $0 \le m < \infty$.
So
$$
\pi_* BP\langle m{-}1 \rangle = BP_*/(v_n \mid n \ge m) \,,
$$
where $BP_* = \Z_{(p)}[v_n \mid n\ge1]$ and $v_0 = p$, and
$$
H_*(BP\langle m{-}1\rangle; \F_p) = \begin{cases}
P(\bar\xi_1^2, \dots, \bar\xi_m^2, \bar\xi_{m+1}, \dots)
	& \text{for $p=2$,} \\
P(\bar\xi_k \mid k\ge1) \otimes E(\bar\tau_k \mid k\ge m)
	& \text{for $p$ odd.}
\end{cases}
$$
The latter is a sub-algebra of the dual Steenrod algebra $A_* = H_*(H\F_p;
\F_p)$.

Suppose that $p$ and $m$ are such that $BP\langle m{-}1\rangle$ admits the
structure of a commutative $S$-algebra.  This is so at least for $m \in
\{0, 1, 2\}$, when $BP\langle-1\rangle = H\F_p$, $BP\langle0\rangle =
H\Z_{(p)}$ and $BP\langle1\rangle = \ell$, respectively, where $\ell$
is the Adams summand of $p$-local connective topological $K$-theory
$ku_{(p)}$.  When $p=2$, $\ell = ku_{(2)}$.

Then the B{\"o}kstedt spectral sequence
$$
E^2_{**} = HH_*(H_*(B; \F_p)) \Longrightarrow H_*(\THH(B); \F_p)
$$
has $E^2$-term
$$
E^2_{**} = \begin{cases}
H_*(BP\langle m{-}1\rangle; \F_2) \otimes E(\sigma\bar\xi_1^2, \dots,
\sigma\bar\xi_m^2, \sigma\bar\xi_{m+1}, \dots)
& \text{for $p=2$,} \\
H_*(BP\langle m{-}1\rangle; \F_p) \otimes E(\sigma\bar\xi_k \mid k\ge1)
\otimes \Gamma(\sigma\bar\tau_k \mid k \ge m)
& \text{for $p$ odd.}
\end{cases}
$$
For $x \in H_*(B; \F_p)$, $\sigma x \in HH_1(H_*(B; \F_p))$ is represented
by the Hochschild $1$-cycle $1 \otimes x$.  The operator $\sigma$ is
a differential ($\sigma^2 = 0$) and a graded derivation ($\sigma(xy) =
x \sigma(y) + (-1)^{|y|} \sigma(x) y$).  Here $\Gamma(-)$ denotes the
divided power algebra.

For $p$ odd, B{\"o}kstedt found differentials
$$
d^{p-1}(\gamma_j(\sigma\bar\tau_k)) = \sigma\bar\xi_{k+1} \cdot
\gamma_{j-p}(\sigma\bar\tau_k)
$$
for $j\ge p$, and in all cases the spectral sequence collapses at the
$E^p$-term.  So
$$
E^\infty_{**} = \begin{cases}
H_*(BP\langle m{-}1\rangle; \F_2) \otimes E(\sigma\bar\xi_1^2, \dots,
\sigma\bar\xi_m^2, \sigma\bar\xi_{m+1}, \dots)
& \text{for $p=2$,} \\
H_*(BP\langle m{-}1\rangle; \F_p) \otimes E(\sigma\bar\xi_1, \dots,
\sigma\bar\xi_m) \otimes P_p(\sigma\bar\tau_k \mid k \ge m)
& \text{for $p$ odd.}
\end{cases}
$$
Here $P_p(-)$ denotes the truncated polynomial algebra of height~$p$.

If $BP\langle m{-}1\rangle$, and thus $\THH(BP\langle m{-}1\rangle)$,
is a commutative $S$-algebra, then there are multiplicative
extensions $(\sigma\bar\xi_k)^2 = \sigma\bar\xi_{k+1}$ for $p=2$ and
$(\sigma\bar\tau_k)^p = \sigma\bar\tau_{k+1}$ for $p$ odd, so
\begin{multline}
H_*(\THH(BP\langle m{-}1\rangle); \F_p) \\
= \begin{cases}
H_*(BP\langle m{-}1\rangle; \F_2) \otimes E(\sigma\bar\xi_1^2, \dots,
\sigma\bar\xi_m^2) \otimes P(\sigma\bar\xi_{m+1})
& \text{for $p=2$,} \\
H_*(BP\langle m{-}1\rangle; \F_p) \otimes E(\sigma\bar\xi_1, \dots,
\sigma\bar\xi_m) \otimes P(\sigma\bar\tau_m)
& \text{for $p$ odd.}
\end{cases}
\end{multline}
For more references and details on the calculation up to this point,
see \cite[Ch.~5]{AnR}.

We now consider the homological homotopy fixed point spectral sequence
for $R = \THH(B)$.  It starts with
$$
E^2_{**} = P(y) \otimes H_*(\THH(B); \F_p)
$$
and by Lemma~\ref{l3.1} it has first differentials
$$
d^2(x) = y \cdot \sigma x
$$
for all $x \in H_*(\THH(B); \F_p)$.  Here $\sigma x \in H_{t+1}(\THH(B);
\F_p)$ is the image of $x \otimes s_1 \in H_t(\THH(B); \F_p) \otimes
H_1(\T; \F_p)$ under the circle action map
$$
\alpha \colon \THH(B) \wedge \T_+ \to \THH(B) \,,
$$
where $s_1 \in H_1(\T; \F_p)$ is the canonical generator.
By Lemma~\ref{l4.6}(a) we have similar differentials $d^2(y^n \cdot x) = y^{n+1}
\cdot \sigma x$ for all $n\ge0$.

Hence we can find the columns of $E^4_{**}$ in the homological homotopy
fixed point spectral sequence by passing to the homology of $E^2_{0,*}
= H_*(\THH(B); \F_p)$ with respect to the operator $\sigma$, at least to
the left of the vertical axis.

\begin{prop}
\label{p6.1}
The homological homotopy fixed point spectral sequence
for $R = \THH(B)$ with $B = BP\langle m{-}1\rangle$, for $p$ and
$m$ such that $B$ is a commutative $S$-algebra, collapses after
the $d^2$-differentials, with the following $E^\infty$-term:

{\rm(a)}\qua
For $p=2$,
$$
E^\infty_{**} = P(y) \otimes P(\bar\xi_1^4, \dots, \bar\xi_m^4,
\bar\xi_{m+1}^2, \xi'_{m+2}, \dots) \otimes E(\bar\xi_1^2
\sigma\bar\xi_1^2, \dots, \bar\xi_m^2 \sigma\bar\xi_m^2)
$$
plus some classes in filtration $s=0$, where $\xi'_{k+1} = \bar\xi_{k+1}
+ \bar\xi_k \sigma\bar\xi_k$ for $k\ge m+1$.

{\rm(b)}\qua
For $p$ odd,
\begin{multline}
E^\infty_{**} = P(y) \otimes
P(\bar\xi_k^p \mid 1 \le k \le m) \otimes P(\bar\xi_{k+1} \mid k \ge m) \\
\otimes E(\tau'_{k+1} \mid k\ge m)
\otimes E(\bar\xi_k^{p-1} \sigma\bar\xi_k \mid 1 \le k \le m)
\end{multline}
plus some classes in filtration $s=0$, where $\tau'_{k+1} = \bar\tau_{k+1}
- \bar\tau_k (\sigma\bar\tau_k)^{p-1}$ for $k \ge m$.
\end{prop}

\begin{proof}
{\rm(a)}\qua
For $B = BP\langle m{-}1\rangle$ and $p=2$ we have
\begin{multline}
E^2_{0,*} = H_*(\THH(BP\langle m{-}1\rangle); \F_2) \\
= P(\bar\xi_1^2, \dots, \bar\xi_m^2, \bar\xi_{m+1}, \dots)
\otimes E(\sigma\bar\xi_1^2, \dots, \sigma\bar\xi_m^2) \otimes
P(\sigma\bar\xi_{m+1}) \,.
\end{multline}
Here $\sigma \colon \bar\xi_k^2 \mapsto \sigma\bar\xi_k^2$ for $1 \le k \le
m$ and $\sigma \colon \bar\xi_{k+1} \mapsto \sigma\bar\xi_{k+1}$ for $k \ge
m$.  We have $\sigma\bar\xi_{k+1} = (\sigma\bar\xi_k)^2$ for $k\ge m+1$.

So the squares $(\bar\xi_k^2)^2 = \bar\xi_k^4$ and $\bar\xi_{m+1}^2$,
as well as the companion classes defined by
$$
\xi'_{k+1} = \bar\xi_{k+1} + \bar\xi_k \sigma\bar\xi_k
$$
for $k \ge m+1$, are $d^2$-cycles, while $E(\bar\xi_k^2,
\sigma\bar\xi_k^2)$ has homology $E(\bar\xi_k^2 \sigma\bar\xi_k^2)$
for each~$k$, and $E(\bar\xi_{m+1}) \otimes P(\sigma\bar\xi_{m+1})$
has homology $\F_2$.

Hence the homological spectral sequence has
$$
E^4_{**} = P(y)
\otimes
P(\bar\xi_1^4, \dots, \bar\xi_m^4, \bar\xi_{m+1}^2, \xi'_{m+2}, \dots)
\otimes
E(\bar\xi_1^2 \sigma\bar\xi_1^2, \dots, \bar\xi_m^2 \sigma \bar\xi_m^2)
$$
plus the image of $\sigma$ in filtration $s=0$.

By our Theorem~\ref{t5.1}(a) applied to the classes $x = \bar\xi_k^2$ for $1
\le k \le m$, in even degree $t = |x|$, the classes $x^2 = \bar\xi_k^4$
and $Q^{t+1}(x) + x \sigma x = \bar\xi_k^2 \sigma\bar\xi_k^2$ are
infinite cycles,  for $Q^{t+1}(\bar\xi_k^2) = 0$ by the Cartan formula.

Similarly, by Theorem~\ref{t5.1}(a) applied to the classes $x = \bar\xi_k$
for $k\ge m+1$, in odd degree $t = |x|$, the classes $x^2 = \bar\xi_k^2$
and $Q^{t+1}(x) + x \sigma x = \bar\xi_{k+1} + \bar\xi_k \sigma\bar\xi_k =
\xi'_{k+1}$ are infinite cycles.  For $Q^{t+1}(\bar\xi_k) = \bar\xi_{k+1}$
by \cite[III.2.2~and~I.3.6]{BMMS86}.

The extra classes in filtration $s=0$ are $y$-torsion, hence infinite
cycles.  Therefore the $E^4$-term above is generated as an algebra
by infinite cycles, so the homological spectral sequence collapses at
this stage.

{\rm(b)}\qua
For $B = BP\langle m{-}1\rangle$ and $p$ odd we have
\begin{multline}
E^2_{0,*} = H_*(\THH(BP\langle m{-}1\rangle); \F_p) \\ = P(\bar\xi_k \mid
k\ge1) \otimes E(\bar\tau_k \mid k\ge m) \otimes E(\sigma\bar\xi_k \mid
1 \le k \le m) \otimes P(\sigma\bar\tau_m) \,.
\end{multline}
Here $\sigma \colon \bar\xi_k \mapsto \sigma\bar\xi_k$ for $1 \le k \le
m$, $\sigma \colon \bar\xi_{k+1} \mapsto 0$ for $k \ge m$ and $\sigma
\colon \bar\tau_k \mapsto \sigma\bar\tau_k$ for $k \ge m$.  We have
$\sigma\bar\tau_{k+1} = (\sigma\bar\tau_k)^p$ for $k \ge m$.

So the $p$-th powers $\bar\xi_k^p$ for $1 \le k \le m$, the classes
$\bar\xi_{k+1}$ for $k\ge m$, and the companion classes defined by
$$
\tau'_{k+1} = \bar\tau_{k+1} - \bar\tau_k (\sigma\bar\tau_k)^{p-1}
$$
for $k \ge m$, are $d^2$-cycles, while $P_p(\bar\xi_k) \otimes
E(\sigma\bar\xi_k)$ has homology $E(\bar\xi_k^{p-1} \sigma\bar\xi_k)$
for each $1 \le k \le m$, and $E(\bar\tau_m) \otimes P(\sigma\bar\tau_m)$
has homology $\F_p$.

Hence the homological spectral sequence has
\begin{multline}
E^4_{**} = P(y) \otimes
P(\bar\xi_k^p \mid 1 \le k \le m) \otimes P(\bar\xi_{k+1} \mid k \ge m) \\
\otimes E(\tau'_{k+1} \mid k\ge m)
\otimes E(\bar\xi_k^{p-1} \sigma\bar\xi_k \mid 1 \le k \le m)
\end{multline}
plus some classes in filtration $s=0$.

Applying our Theorem~\ref{t5.1}(b) to the classes $x = \bar\xi_k$ for $1
\le k \le m$, in even degree $t = |x|$, the classes $x^p = \bar\xi_k^p$
and $x^{p-1} \sigma x = \bar\xi_k^{p-1} \sigma\bar\xi_k$ are
infinite cycles.

Similarly, applying Theorem~\ref{t5.1}(c) to the classes $x = \bar\tau_k$ for
$k\ge m$, in odd degree $t = |x| = 2p^k-1$, the classes $\beta Q^{p^k}(x)
= \bar\xi_{k+1}$ and $Q^{p^k}(x) - x (\sigma x)^{p-1} = \bar\tau_{k+1} -
\bar\tau_k (\sigma\bar\tau_k)^{p-1} = \tau'_{k+1}$ are infinite cycles,
for $Q^{p^k}(\bar\tau_k) = \bar\tau_{k+1}$ and $\beta \bar\tau_{k+1} =
\bar\xi_{k+1}$ by \cite[III.2.3~and~I.3.6]{BMMS86}.

Hence the $E^4$-term above is generated as an algebra by infinite
cycles, and the homological spectral sequence collapses after the
$d^2$-differentials.
\end{proof}

For convenience in the comparison with $ko$, we make the case $B = ku$
at $p=2$ explicit:

\begin{cor}
The homological homotopy fixed point spectral sequence for $R
= \THH(ku)$ at $p=2$ collapses after the $d^2$-differentials,
with
$$
E^\infty_{**} = P(y) \otimes P(\bar\xi_1^4, \bar\xi_2^4, \bar\xi_3^2,
\xi'_4, \dots) \otimes E(\bar\xi_1^2 \sigma\bar\xi_1^2, \bar\xi_2^2
\sigma\bar\xi_2^2)
$$
plus some classes in filtration $s=0$, where $\xi'_{k+1} = \bar\xi_{k+1}
+ \bar\xi_k\sigma\bar\xi_k$ for $k\ge3$.
\end{cor}

\begin{prop}
The homological homotopy fixed point spectral sequence for $R = \THH(B)$
collapses after the $d^2$-differentials, in both of the cases:

{\rm(a)}\qua
$B = ko$ and $p=2$, when
$$
E^\infty_{**} = P(y) \otimes P(\bar\xi_1^8, \bar\xi_2^4, \bar\xi_3^2,
\xi'_4, \dots) \otimes E(\bar\xi_1^4 \sigma\bar\xi_1^4, \bar\xi_2^2
\sigma\bar\xi_2^2)
$$
plus classes on the vertical axis,
and

{\rm(b)}\qua
$B = \tmf$ and $p=2$, when
$$
E^\infty_{**} = P(y) \otimes P(\bar\xi_1^{16}, \bar\xi_2^8, \bar\xi_3^4,
\bar\xi_4^2, \xi'_5, \dots) \otimes E(\bar\xi_1^8 \sigma\bar\xi_1^8,
\bar\xi_2^4 \sigma\bar\xi_2^4, \bar\xi_3^2 \sigma\bar\xi_3^2)
$$
plus classes on the vertical axis.
\end{prop}

\begin{proof}
{\rm(a)}\qua
For $B = ko$ with $H_*(B; \F_2) = (A/\!/A_1)_* =
P(\bar\xi_1^4, \bar\xi_2^2, \bar\xi_3, \dots)$ we have
$$
H_*(\THH(ko); \F_2) = P(\bar\xi_1^4, \bar\xi_2^2, \bar\xi_3, \dots)
\otimes E(\sigma\bar\xi_1^4, \sigma\bar\xi_2^2) \otimes
P(\sigma\bar\xi_3) \,.
$$
See \cite[6.2(a)]{AnR}.

As in the proof of Proposition~\ref{p6.1}, the squares $\bar\xi_1^8$,
$\bar\xi_2^4$ and $\bar\xi_3^2$, as well as the classes $\xi'_{k+1} =
\bar\xi_{k+1} + \bar\xi_k \sigma\bar\xi_k$ for $k\ge3$ are $d^2$-cycles,
while $E(\bar\xi_1^4, \sigma\bar\xi_1^4)$ and $E(\bar\xi_2^2,
\sigma\bar\xi_2^2)$ have homology $E(\bar\xi_1^4 \sigma\bar\xi_1^4)$
and $E(\bar\xi_2^2 \sigma\bar\xi_2^2)$, respectively.  The homology of
$E(\bar\xi_3) \otimes P(\sigma\bar\xi_3)$ is $\F_2$.  So
$$
E^4_{**} = P(y) \otimes P(\bar\xi_1^8, \bar\xi_2^4, \bar\xi_3^2,
\xi'_4, \dots) \otimes E(\bar\xi_1^4 \sigma\bar\xi_1^4, \bar\xi_2^2
\sigma\bar\xi_2^2)
$$
plus some classes in filtration $s=0$.

By Theorem~\ref{t5.1}(a), all of these algebra generators are in fact
infinite cycles, so the homological spectral sequence collapses,
as claimed.

{\rm(b)}\qua
For $B = \tmf$ with $H_*(B; \F_2) = (A/\!/A_2)_* =
P(\bar\xi_1^8, \bar\xi_2^4, \bar\xi_3^2, \bar\xi_4, \dots)$ we have
$$
H_*(\THH(\tmf); \F_2) = P(\bar\xi_1^8, \bar\xi_2^4, \bar\xi_3^2, \bar\xi_4,
\dots) \otimes E(\sigma\bar\xi_1^8, \sigma\bar\xi_2^4, \sigma\bar\xi_3^2)
\otimes P(\sigma\bar\xi_4) \,.
$$
See \cite[6.2(b)]{AnR}.  This gives the $E^2$-term of the homological
spectral sequence, and as before its homology with respect to the
$\sigma$-operator is
$$
E^4_{**} = P(y) \otimes P(\bar\xi_1^{16}, \bar\xi_2^8, \bar\xi_3^4,
\bar\xi_4^2, \xi'_5, \dots) \otimes E(\bar\xi_1^8 \sigma\bar\xi_1^8,
\bar\xi_2^4 \sigma\bar\xi_2^4, \bar\xi_3^2 \sigma\bar\xi_3^2)
$$
plus some classes in filtration $s=0$.

By Theorem~\ref{t5.1}(a), all of these algebra generators are in fact
infinite cycles, so the homological spectral sequence collapses,
as claimed.
\end{proof}

\begin{thm}
\label{t6.4}
The homological homotopy fixed point spectral sequence for $R = \THH(B)$
collapses after the $d^2$-differentials, in both of the cases:

{\rm(a)}\qua
$B = MU$, with
$$
E^\infty_{**} = P(y) \otimes P(b_k^p \mid k\ge1) \otimes
E(b_k^{p-1} \sigma b_k \mid k\ge1)
$$
plus classes in filtration zero, and

{\rm(b)}\qua
$B = BP$, with
$$
E^\infty_{**} = P(y) \otimes P(\xi_k^p \mid k\ge1) \otimes
E(\xi_k^{p-1} \sigma\xi_k \mid k\ge1)
$$
plus classes in filtration zero.  (When $p=2$, substitute $\xi_k^2$
for $\xi_k$.)
\end{thm}

Note that we do not need to assume that $BP$ is a commutative
$S$-algebra for the result in part~(b).

\begin{proof}(a)\qua
The integral homology algebra of $MU$ is $H_*(MU; \Z) = \Z[b_k \mid
k\ge1]$, where $b_k$ in degree $2k$ is the stabilized image of the
generator $\beta_{k+1} \in H_{2k+2}(BU(1); \Z)$, under the zero-section
identification $BU(1) \simeq MU(1)$.  So
$$
H_*(MU; \F_p) = P(b_k \mid k\ge1)
$$
is concentrated in even degrees,
and the $E^2$-term of the B{\"o}kstedt spectral sequence is
$$
E^2_{**} = HH_*(H_*(MU; \F_p)) = H_*(MU; \F_p) \otimes
E(\sigma b_k \mid k\ge1) \,.
$$
All the algebra generators are in filtrations $s \le 1$, so the spectral
sequence collapses at this stage.  There are no algebra extensions, since
for $p=2$, $(\sigma b_k)^2 = Q^{2k+1}(\sigma b_k) = \sigma Q^{2k+1}(b_k)
= 0$, where $Q^{2k+1}(b_k) = 0$ because it has odd degree.  For $p$
odd, $(\sigma b_k)^2 = 0$ by graded commutativity, because $\sigma b_k$
has odd degree.  Thus
$$
H_*(\THH(MU); \F_p) = H_*(MU; \F_p) \otimes E(\sigma b_k \mid k\ge1)
\,.
$$
This much can also be read off from \cite[4.3]{MS93}, or from Cohen and
Schlichtkrull's formula $\THH(MU) \simeq MU \wedge SU_+$ \cite{CS}.

The homological homotopy fixed point spectral sequence has $E^2$-term
$$
E^2_{**} = P(y) \otimes P(b_k \mid k\ge1) \otimes E(\sigma b_k
\mid k\ge1) \,.
$$
Its homology with respect to the $d^2$-differential, satisfying
$d^2(b_k) = y \cdot \sigma b_k$, is
$$
E^4_{**} = P(y) \otimes P(b_k^p \mid k\ge1) \otimes
E(b_k^{p-1} \sigma b_k \mid k\ge1)
$$
plus the usual $y$-torsion on the vertical axis.  By Theorem~\ref{t5.1}(a)
and~(b), the algebra generators of this $E^4$-term are all infinite
cycles.  Hence the spectral sequence collapses at this stage.

{\rm(b)}\qua
The Brown--Peterson spectrum $BP$ was originally constructed to have
mod~$p$ homology
$$
H_*(BP; \F_p) = \begin{cases}
P(\xi_k^2 \mid k\ge1) & \text{for $p=2$,} \\
P(\xi_k \mid k\ge1) & \text{for $p$ odd.}
\end{cases}
$$
This equals the sub-algebra $(A/\!/E)_*$ of $A_*$ that is dual to the
quotient algebra $A/\!/E = A/A\beta A$ of $A$.  Hereafter we focus on
the odd-primary case; the reader should substitute $\xi_k^2$ for $\xi_k$
when $p=2$.

The spectrum $BP$ is known to be an (associative) $S$-algebra, and to
receive an $S$-algebra map from $MU$ \cite[3.5]{BJ02}.  This map induces
a split surjective algebra homomorphism $H_*(MU; \F_p) \to H_*(BP;
\F_p)$ in homology, which maps $b_{p^k-1}$ to $\xi_k$ for $k\ge1$
and takes the remaining algebra generators $b_i$ to $0$ for $i \ne
p^k-1$.  For the homology of $BP$ injects into $H_*(H\Z_{(p)}; \F_p)$
and at the level of second spaces the composite map of spectra $MU \to
BP \to H\Z_{(p)}$ is a $p$-local equivalence $MU(1) \to K(\Z_{(p)},
2)$.  The generator $\beta_{i+1} \in \tilde H_{2i+2}(MU(1); \F_p)$
maps to $b_i \in H_{2i}(MU; \F_p)$, while the corresponding generator
$\beta_{i+1} \in \tilde H_{2i+2}(K(\Z_{(p)}, 2); \F_p)$ maps to $\xi_k
\in H_{2i}(H\Z_{(p)}; \F_p)$ when $i = p^k-1$ and to $0$ otherwise
\cite[Ch.~5]{Mi58}.  This proves the claim.

The B{\"o}kstedt spectral sequence for $BP$ has $E^2$-term
$$
E^2_{**} = HH_*(H_*(BP; \F_p)) = H_*(BP; \F_p)
\otimes E(\sigma\xi_k \mid k\ge1) \,.
$$
Note that the map $MU \to BP$ induces a surjection of B{\"o}kstedt
spectral sequence $E^2$-terms.  Thus the fact that the B{\"o}kstedt
spectral sequence for $MU$ collapses at $E^2$ with no algebra extensions
implies the corresponding statement for $BP$, also without the assumption
that $BP$ is a commutative $S$-algebra.  We can conclude that
$$
H_*(\THH(BP); \F_p) = H_*(BP; \F_p) \otimes E(\sigma\xi_k \mid k\ge1)
\,.
$$

The homological homotopy fixed point spectral sequence has $E^2$-term
$$
E^2_{**} = P(y) \otimes P(\xi_k \mid k\ge1)
\otimes E(\sigma\xi_k \mid k\ge1) \,.
$$
Again the map $MU \to BP$ induces a surjection of $E^2$-terms, so the
$d^2$-differentials satisfy $d^2(\xi_k) = y \cdot \sigma\xi_k$ and $d^2(y)
= 0$, and are derivations.  This leaves
$$
E^4_{**} = P(y) \otimes P(\xi_k^p \mid k\ge1)
\otimes E(\xi_k^{p-1} \sigma\xi_k \mid k\ge1)
$$
plus some $y$-torsion on the vertical axis, and the map from the
$E^4$-term of the spectral sequence for $MU$ is still surjective.
Thus the spectral sequence for $BP$ also collapses at this stage.
\end{proof}

\section{ Generalizations and comments }

In this section we note some generalizations of our results, and also
comment on the relation to similar patterns of differentials in other
spectral sequences.  The generalizations are of two sorts.  First, we can
replace the homotopy fixed points construction by the Tate construction
or the homotopy orbits.  Second, we can change the group of equivariance.
We consider these in order.

First, there are spectral sequences similar to the one considered here
for the Tate construction $X^{t\T} = [\widetilde{E\T} \wedge F(E\T_+,
X)]^{\T}$ (denoted $t_{\T}(X)^{\T}$ in \cite{GM95} and $\widehat{\mathbb
H}(\T, X)$ in \cite{AuR02}) and the homotopy orbit spectrum $X_{h\T} =
E\T_+ \wedge_{\T} X$.

\begin{prop}
There is a natural spectral sequence
$$
\widehat{E}^2_{**} = 
\widehat{H}^{-*}(\T; H_*(X; \F_p)) = P(y,y^{-1}) \otimes H_*(X; \F_p)
$$
with $y$ in bidegree $(-2,0)$, which is conditionally convergent to
the continuous homology $H^c_*(X^{t\T}; \F_p)$.  We call this the {\sl
homological Tate spectral sequence}.  If $H_*(X; \F_p)$ is bounded below
and finite in each degree, or the spectral sequence collapses at a finite
stage, then the spectral sequence is strongly convergent.
\end{prop}

\begin{prop}
\label{p7.2}
There is a natural spectral sequence
$$
E^2_{**} = H_{*}(\T; H_*(X; \F_p)) = P(y^{-1}) \otimes H_*(X; \F_p)
$$
with $y^{-1}$ in bidegree $(2,0)$, which is strongly convergent to
$H_*(X_{h\T}; \F_p)$.  We call this the {\sl homological homotopy orbit
spectral sequence}.  (Note that for $X_{h\T}$ the continuous homology
is the same as the ordinary homology.)
\end{prop}

Further, the middle and right hand maps of the (homotopy) norm cofiber
sequence
$$
\Sigma X_{h\T} \xrightarrow{N} X^{h\T} \to X^{t\T} \to \Sigma^2 X_{h\T}
$$
induce the homomorphisms of $E^2$-terms given by tensoring $H_*(X; \F_p)$
with the short exact sequence of $P(y)$-modules
$$
0 \to P(y) \to P(y,y^{-1}) \to \Sigma^2 P(y^{-1}) \to 0.
$$
Thus the homological Tate spectral sequence is a full-plane spectral
sequence whose $E^2$-term is obtained by continuing the $y$-periodicity
in the homological homotopy fixed point spectral sequence into the
right half-plane, and the homological homotopy orbit spectral sequence
(shifted~$2$ degrees to the right from Proposition~\ref{p7.2}) has the quotient
of these as its $E^2$-term.

Proposition~\ref{p4.4} and Theorem~\ref{t5.1} apply equally well to all three spectral
sequences.  For details, see the thesis of Lun{\o}e--Nielsen \cite{L-N}.

Second, we could also consider these three spectral sequences for the
action of a finite cyclic subgroup $C$ of $\T$.  For example, there is
the homological Tate spectral sequence
$$
\widehat{E}^2_{**} = \widehat{H}^{-*}(C; H_*(X; \F_p))
$$
converging conditionally to $H^c_*(X^{tC}; \F_p)$.  The analogue of
Lemma~\ref{l4.6} still holds, so that there are isomorphisms
$$
\widehat{E}^r_{**} \cong \widehat{H}^{-*}(C; \F_p) \otimes
\widehat{E}^r_{0,*}
$$
for all $r\ge2$ (and now $y$ is invertible, so there is no $y$-torsion),
and all differentials are determined by those originating on the
vertical axis $\widehat{E}^r_{0,*}$.  In turn, the latter differentials
are determined by those in the $\T$-equivariant case, by naturality
with respect to the restriction map $X^{t\T} \to X^{tC}$.  Therefore
the collapse results in Theorem~\ref{t5.1} also hold in these cases.
See \cite{L-N} for more details.

These latter spectral sequences, for finite subgroups $C \subset \T$,
are essential in the analysis of the topological model $TF(B)$ for the
negative cyclic homology of~$B$, and the topological cyclic homology
$TC(B)$.

Though the differentials here allow us to determine $E^\infty_{**}$ in
the cases of interest (see Section~6), there are still $A_*$-comodule
extensions hidden by the filtration.  These are of course of critical
importance for the analysis of the Adams spectral sequence~(\ref{e1.3}).
A more elaborate study of the geometry of the universal examples used
in Section~5 allows these to be recovered.  This too can be found in
\cite{L-N}.

Finally, it is interesting to compare the formulas for differentials
here to analogous results in other spectral sequences.  The first to be
considered was the Adams spectral sequence, where the results are due to
Kahn \cite{Ka70}, Milgram \cite{Mi72}, M\"{a}kinen \cite{Mak73} and the
first author \cite[Ch.~VI]{BMMS86}.  For simplicity, let us assume $p=2$
in this discussion, as there are several cases to be considered at odd
primes (\cite[VI.1.1]{BMMS86}).  Suppose that $x$ is in the $E_r$-term
of the Adams spectral sequence
$$
E_2^{**} = \Ext_A^{**}(H^*(R; \F_2), \F_2)
	\Longrightarrow \pi_*(R)^\wedge_2 \,,
$$
where $R$ is a commutative $S$-algebra.  The commutative $S$-algebra
structure of $R$ induces Steenrod operations in the $E_2$-term of the
Adams spectral sequence, which are the analog in this situation of the
Dyer--Lashof operations in $H_*(R; \F_2)$.  (In fact, under the Hurewicz
homomorphism, they map to the Dyer--Lashof operations.)  Then, in most
cases we have
\begin{equation}
d_*(Sq^j x) = Sq^j d_r(x) \,\,\, {\dot{+}} \,\,\, a Sq^{j-v} x \,,
\label{e7.3}
\end{equation}
where $A\,\,\dot{+}\,\,B$ denotes whichever of $A$ or $B$ is in the
lower filtration, or their sum, if they are in the same filtration.
The subscript in $d_*$ is then the difference in filtrations between the
right and left hand sides.  In this formula, $a$ is an infinite cycle in
the Adams spectral sequence for the homotopy groups of spheres, and $a$
and $v$ are determined by $j$ and the degree of $x$.  When the first
half of the right hand side dominates we have
$$
d_{2r-1}(Sq^j x) = Sq^j d_r(x) \,,
$$
and this formula resembles the formula 
$$
d^{2r}(\beta^\epsilon Q^i(x)) = \beta^\epsilon Q^i(d^{2r}(x))
$$
of Proposition~\ref{p4.4}, in that both essentially say that the relevant
differential commutes with the  Dyer--Lashof operations.  The fact that
the length of the differential increases  from $r$ to $(2r-1)$ when we
apply the squaring operation in the Adams spectral sequence reflects
the difference between the homotopy fixed point filtration and the Adams
filtration, and the way in which they interact with the extended powers.
A more extreme difference occurs when the second term $a Sq^{j-v} x$
is involved.  In the homological homotopy fixed point spectral sequence
this term disappears, essentially because the element $a \in \pi_* S$
is mapped to $0$ by the Hurewicz homomorphism.  Homotopical homotopy
fixed point spectral sequences, as in \cite{AuR02}, will have differential
formulas with two parts, as in the Adams spectral sequence.  Such two part
formulas for differentials reflect {\it universally hidden extensions\/}
in the following sense.

The differential~(\ref{e7.3}) arises from decomposing the boundary of the
cell on which $Sq^j x$ is defined into two pieces.  One of the pieces
carries $Sq^j d_r(x)$ and the other carries $a Sq^{j-v} x$.  The half
that lies in the lower filtration is killed by the differential (\ref{e7.3}),
and therefore appears to be~$0$ in the associated graded $E_\infty$-term.
However, the geometry of the situation shows that it is actually equal
to the half of the formula that lies in the higher filtration, modulo
still higher filtrations.  Thus we have a universally hidden extension,
that is, an expression which is 0 in the associated graded, by virtue
of being equal to an expression which lies in a higher filtration.
We should expect this sort of phenomenon to occur in homotopical homotopy
fixed point spectral sequences.

Finally, Theorem~\ref{t5.1} seems to be particular to the homological
homotopy fixed point spectral sequence.  Certainly the Adams spectral
sequence seems to have no analog of this extreme cutoff, in which certain
terms die at $E^r$ and the remaining terms live to $E^\infty$.

\makeatletter
\def\thebibliography#1 {\@thebibliography@{AuR02}\small\parskip0pt 
plus2pt\relax\itemsep 1pt plus 1pt
\addcontentsline{toc}{section}{Bibliography}}
\makeatother

\Addresses\recd

\end{document}